\documentclass[a4paper,12pt,oneside,reqno]{amsart}
\usepackage[latin1]{inputenc}
\usepackage[english]{babel}
\usepackage[dvips]{graphicx}
\usepackage{amssymb}
\usepackage{amsmath}
\usepackage{amsthm}
\usepackage{url}
\usepackage{hyperref}
\usepackage{a4wide}

\usepackage[hmargin=2.2cm,vmargin=2.3cm]{geometry} 

\def\[#1\]{\begin{align*}#1\end{align*}}
\def\be#1\ee{\begin{align}#1\end{align}}
\def\bea#1\eea{\begin{align}#1\end{align}}
\def\ben#1\een{\begin{align*}#1\end{align*}}

\newcommand{\n}[1]{\left\Vert #1\right\Vert}
\newcommand{\la}{\left\langle}
\newcommand{\ra}{\right\rangle}
\newcommand{\R}{\mathbb{R}}

\newcommand{\N}{\mathbb{N}}

\newcommand{\mc}[1]{\mathcal{#1}}
\newcommand{\s}{\mc}
\newcommand{\p}[2]{\frac{\partial #1}{\partial #2}}

\def\ol#1{\overline{#1}}

\def\[#1\]{\begin{align*}#1\end{align*}}
\def\be#1\ee{\begin{align}#1\end{align}}
\def\bea#1\eea{\begin{align}#1\end{align}}
\def\ben#1\een{\begin{align*}#1\end{align*}}

\newtheoremstyle{theorem}{0.5cm}{0.5cm}%
   {}
   {}
   {\bfseries}
   {}
   {2ex}
   {\thmname{#1}\thmnumber{ #2}\thmnote{ #3}}
\theoremstyle{theorem}
\newtheorem{theorem}{Theorem}[section]

\newtheorem{example}[theorem]{Example}

\newtheorem{remark}[theorem]{Remark}

\newtheorem{lemma}[theorem]{Lemma}

\begin{document}
\title[System of transport equations and  inverse problem]
{A system of transport equations and an  inverse problem from radiation therapy}

\author{Jouko Tervo}
\address{University of Eastern Finland, Department of Applied Physics, P.O.Box 1627, 70211 Kuopio, Finland}
\email{040tervo@gmail.com}
 
\date{\today}

\maketitle

\begin{abstract}

The paper considers existence results of solution for a linear coupled system of  Boltzmann transport equations and related inverse problem. The system models the evolution of three species of particles, photons, electrons and positrons.  Hyper-singularities of differential cross sections associated with charged particle transport cause that  modelling contains the first order partial differential term with respect to energy and the second order partial differential term with respect to velocity angle. 
The overall system is a partial integro-differential equation.
The model is intended especially for dose calculation in the forward problem  of radiation therapy. Firstly we consider a single transport equation for charged particles. After that
we verify under physically relevant assumptions that the coupled equation together with relevant initial and inflow boundary values has a unique solution in appropriate $L^2$-based spaces. The existence of solutions for the adjoint problem is verified as well.
Moreover, we deal with the related inverse problem, a so called inverse radiation treatment planning problem. It is formulated as an optimal boundary control problem.  Variational equations for an optimal control related to an appropriate differentiable strictly convex object function are verified. Its solution can be used as an initial point for an actual  optimization which needs global optimization methods.

\end{abstract}


{\small \textbf {AMS-Classification:}  35Q20, 45E99, 35R30, 49N45} 

\maketitle

\section{Introduction}\label{intro}

The Boltzmann transport equation (BTE) is an integro-partial differential equation which physically is based on the conservation laws. It has applications in many fields of scientific computation, including in among others optical tomography, cosmic radiation,
nanotechnology (e.g. plasma physics) and radiation therapy, which is considered in this paper.
For general mathematical theory of BTE with relevant boundary conditions we refer to \cite{agoshkov} and \cite{dautraylionsv6}. See also \cite{casezweifel}, \cite{cercignani}, \cite{duderstadt79},  
\cite{pomraning} where the subject is considered from more physical point of view.
Some more recent issues related to BTE can be found \cite{mokhtarkharroubi}, and some non-linear aspects in \cite{bellamo}, \cite{ukai}.  A  mathematical survey of non-linear collision theory of particle physics (especially in dilute gases)  is found in \cite{villani}. 
Finally, for topics related to Monte Carlo methods in the context of BTE,
both from theoretical and practical points of view,
we refer to \cite{lapeyre03} and \cite{seco13}.

The goal in radiation therapy is to
generate dose distributions in such a way that the prescribed dose conforms
as well as possible to the target volume,
while healthy tissue, and especially
the so-called critical organs, achieve as low dose as possible.
In {\it external radiotherapy} the problem is to find the optimal dose by defining the field intensity, that is the incoming particle flux, on (patches of) the patient surface. 
In {\it internal radiotherapy} the radioactive sources are to be
positioned inside the patient tissue such that the desired dose distribution is achieved.
The determination of the external particle flux (or the distribution of internal sources) required to deliver the desired dose distribution is called
the \emph{inverse radiation treatment planning} (IRTP)  problem,
which from mathematical point of view is an {\it inverse problem}.
The calculation of particle fluxes or dose in the patient tissue when the incoming fluxes or internal sources are known, is called {\it dose calculation},
and it is considered as a {\it forward problem}.

The solution of IRTP always requires some \emph{dose calculation model} which in this paper is founded on the relevant system of BTEs.
In radiation therapy the BTE-system describes how radiation is scattered and absorbed in a 
tissue. High energy incoming particles, such as photons or electrons  mobilize
(at least) three kinds of particles, photons, electrons 
and positrons, whose simultaneous evolution should be taken into account in the 
transport model.
The potential creation of (or contamination by)
other heavy particles (such as neutrons)
will not be taken into account since their
contribution is negligible (cf. \cite{seco13}). 
Dose calculation models governed by the  BTE-systems are valid in inhomogeneous material. They take rigorously into account the scattering and absorption effects in physically solid way.

We assume  that the transport of  particles is ruled by the following \emph{linear coupled system of three BTEs}
for $\psi=(\psi_1,\psi_2,\psi_3)$ on $G\times S\times I$,
\bea
&
\omega\cdot\nabla_x\psi_1+\Sigma_1\psi_1-K_{r,1}\psi=f_1, \label{intro1}\\
&
a_j(x,E){\p {\psi_j}{E}}
+b_j(x,E)\Delta_S\psi_j
+\omega\cdot\nabla_x\psi_j+\Sigma_{j}\psi_j-K_{r,j}\psi
=f_j,\quad j=2,3.\label{intro2}
\eea
together with the inflow boundary condition on $\Gamma_-$,
\be\label{intro3}
{\psi_j}_{|\Gamma_-}=g_j,\quad j=1,2,3,
\ee
and initial value  condition on $G\times S$,
\be\label{intro4}
\psi_j(\cdot,\cdot,E_{\rm m})=0,\quad j=2,3.
\ee
Above $K_{r,j}$ are the \emph{restricted collision operators} which are assumed to be of the form
\be\label{intro5}
(K_{r,j}\psi)(x,\omega,E)=\sum_{k=1}^3\int_{I'}\int_{S'}\sigma_{kj}(x,\omega',\omega,E',E)\psi_k(x,\omega',E') d\rho^{kj}_{S}(\omega')d\rho^{kj}_I(E'),\quad j=1,2,3,
\ee
where $d\rho^{kj}_{S}(\omega')$ and $d\rho^{kj}_I(E')$  are relevant pairs of measures corresponding to the different
varieties of restricted collision operators (see section \ref{coup-ex} below).
Here and below we denote by $I'$ and $S'$ the sets $I$ and $S$ when the integration variable is $E'$ and $\omega'$, respectively.
On the right, the functions $f_j$ represent the (internal) sources,
and $g_j$ in \eqref{intro3} are boundary (external) sources.
$G\subset\R^3$ is the spatial (patient) domain, $S\subset\R^3$ the unit sphere
(velocity angle domain) and $I=[E_0,E_m]\subset\R$ the energy interval, where $E_{\rm m}$ is the so called \emph{cut-off energy}. $b_j(x,E)\Delta_S$ is the \emph{Laplace-Beltrami operator} on $S$ together with (angular) diffusion  coefficient.
The functions $\Sigma_j=\Sigma_{j}(x,E)$ are the \emph{restricted total cross-sections}, 
$\sigma_{kj}=\sigma_{kj}(x,\omega',\omega,E',E)$ are the \emph{restricted differential cross-sections},
and the factors $a_j=a_{j}(x,E)$ are the so-called \emph{restricted stopping powers}. 
The  system is \emph{coupled} through the collision operators $K_{r,j}$. 
The solution $\psi=(\psi_1,\psi_2,\psi_3)$ of the problem \eqref{intro1}-\eqref{intro4} is a vector-valued function whose components
describe the particle number densities of photons, electrons and positrons, respectively.
The equations (\ref{intro1}), (\ref{intro2}) are steady state (stationary state) equations where there is no time dependence.
This is basically always the case in applications related to radiotherapy,
because the relevant radiation field $\psi$ anyway reach the steady state
nearly instantly (\cite{borgers99}).
The use of BTEs in dose calculation needs the choice of { total and differential cross-sections}.
In radiation therapy the cross-sections of primary interest are those for water (tissue), bone and air (void-like regions). 
For a more thorough discussion on the cross sections relevant to radiation therapy, we refer to \cite{bomanthesis}, \cite{hensel}.

We point out   that  for some charged particle interactions
the differential cross-sections 
may have singularities or even {\it hyper-singularities} which 
imply that the model equations  (\ref{intro2}) includes
the partial differential terms $a_j(x,E){\p {\psi_j}{E}}$, 
$b_j(x,E)\Delta_S\psi_j$ (\cite{tervo18}, \cite{tervo19}). The equations are approximations of the exact partial hyper-singular integro-differential equations related to M\o ller-type interactions. 
In literature the approximations to cover these problematic interactions are called the  {Continuous Slowing Down-Boltzmann Transport Equations}
(CSDA-BTE).
We notice
that the above  transport model is linear and, therefore, neglects any non-linear interactions which have been reported to be meaningless in dose calculation. 
In addition, the inflow boundary condition (\ref{intro3}) is not exactly correct because minor part of the  particles return to the patient domain $G$.
The reflection boundary conditions of the form
${\psi_j}_{|\Gamma_-}=R_j({\psi_j}_{|\Gamma_+})+g_j,\ j=1,2,3$, where $R_j$ are appropriate (unbounded) reflection operators would be more accurate (see \cite{tervo18-up},
section 7.4,  \cite{dautraylionsv6}, pp. 251-262).

The (classical) example of the initial inflow boundary value  problems of the form  (\ref{intro1})-(\ref{intro4}) is the case where only the equation like (\ref{intro1}) is considered.
Existence of solutions for this kind of  transport problems
has been studied for  single equations e.g. in \cite{egger14,dautraylionsv6,agoshkov} and
for  coupled systems   \cite{tervo17-up}. 
In the references it is assumed that the restricted collision operator $K_r$  satisfies a (partial) \emph{ Schur criterion } for boundedness (\cite{halmos}, p. 22).
In \cite{tervo17,tervo18,frank10} one has  studied 
existence of solutions for the case where the term $b_j(x,E)\Delta_S$ vanishes.
For related existence results of (deterministic and linear) Fokker-Planck type equations we refer to \cite[Appendix A]{degond}, \cite{tian}, \cite{degond},
\cite{le-bris}, \cite{chupin} and
\cite{herty11} (where existence results are exposed in the context of dose calculation  for  radiation treatment planning).  
In part of the references the state space is $G\times V$ where $V$ is velocity space but they can be formulated for $G\times S\times I$ via the (local)
diffeomorphism 
$h(\omega,E):=\sqrt{E}\omega$ (see Appendix below).

In \cite{tervoarxiv-18} we derived a refined expression for the exact transport operator related to a charged particle transport. The analysis therein was carried out only for the M\o ller-type scattering which is a special  prototype of hyper-singular interactions. 
This interaction models the electron's (and positron's) inelastic collisions and other type  of collisions (such as to Bremsstrahlung) could be handled analogously. 
In section \ref{tr-app} we recall from \cite{tervo19} the relevant approximation of the exact hyper-singular transport operator. 
The analysis yields a foundation for the CSDA-Fokker-Planck type approximations (\ref{intro2}) of transport operators.
The restricted collision operator $K_r$ is roughly speaking  the residual when the singular part is separated from the hyper-singular collision operator.
In \cite{tervo18} we showed  that $K_r$ is a {bounded operator} in $L^2(G\times S\times I)$.  In addition, we showed that $\Sigma-K_r$ is \emph{coercive (accretive) operator} in $L^2(G\times S\times I)$ under relevant assumptions.

In section \ref{trunc-exists} we consider the existence of solutions. At first we recall from \cite{tervo19} existence results in $L^2(G\times S\times I)$-based spaces for the transport problem corresponding to a single charged particle transport. 
The formulations
for the above coupled system 
\eqref{intro1}-\eqref{intro4} in $L^2(G\times S\times I)^3$-based spaces are given in section \ref{cosyst}. 
The transport problems  are analysed in the variational (weak) form 
\[
\widetilde B(\psi,v)=Fv,\ v\in \widehat{\mc H}.
\]
Here $\widehat{\mc H}$ is a suitable space of test functions and $\widetilde B(.,.)$ and $F$ are the due bilinear  and linear forms, respectively. 
The variational formulation is an essential step in order  to show existence results e.g. by applying Lions-Lax-Milgram Theorem based methods. It is also needed in searching numerical solutions for the forward problem by \emph{finite element methods} (FEM).

In section \ref{irtpre} we consider the above mentioned IRTP problem under physical criteria for optimization (for some general backgrounds see e.g. \cite{shepard99}, \cite{mayles}).  
The patient domain $G\subset\R^3$ consists of tumor volume ${\s T}$, critical organ's region ${\s C}$ and the normal tissue's region ${\s N}$. 
Hence  $G={\s T}\cup {\s C}\cup {\s N}$
where the union is mutually disjoint.
The tumor volume (that is, the target) includes the tumor and some safety margin. Critical organs and normal tissue are build up of healthy tissue,
and should receive as low a dose as possible.
Typically the resulting object function  based on the physical criteria is of the form (see section \ref{objectf})
\bea\label{io}
J({\bf f},{\bf g})=&c_{\s T}\n{D_0-{\s D}({\bf f},{\bf g})}_{L^2({\s T})}^2
+
c_{\s C}\n{(D_C-{\s D}({\bf f},{\bf g}))_-}_{L^2({\s C})}^2 \nonumber \\
&+
c_{\s N}\n{(D_N-{\s D}({\bf f},{\bf g}))_-}_{L^2({\s N})}^2 \nonumber \\
&+
c_{\rm DV}\Big(\Big(v_C-\frac{1}{m({\s C})}\int_{\s C}H(({\s D}({\bf f},{\bf g}))(x)-d_C) dx\Big)_-\Big)^2
\eea
where ${\bf f}=(f_1,f_2,f_3),\ {\bf g}=(g_1,g_2,g_3)$ and
where $c_{\s T}, c_{\s C}, c_{\s N}, c_{\rm DV}$ are positive weights,
$m$ is the 3-dimensional Lebesgue measure, $H$ is the Heaviside function
and $a_-$ denotes the negative part of $a\in\R$.
Here  ${\s D}({\bf f},{\bf g})=D(\psi({\bf f},{\bf g}))$ where $D$ is the dose operator
(see section \ref{RTP})
and $\psi=\psi({\bf f},{\bf g})$ is the solution of (\ref{intro1})-(\ref{intro4}) with data $({\bf f},{\bf g})$.
We note that ${\bf f}=0$ for external therapy and ${\bf g}=0$ for internal therapy. 
In (\ref{io}) the first three terms are convex and (locally) Lipschitz continuous and the first term is also differentiable. The last term is both non-convex and non-differentiable in general but
it can be replaced a by Lipschitz continuous counterpart by replacing Heaviside function $H$  by a Lipschitz continuous approximation.
After this replacement the whole object function (\ref{io}) is (locally) Lipschitz continuous. In addition, the admissible sets (as given in section \ref{RTP}) are convex.

In practice, solving deterministically  the discretized BTE is a challenging numerical task because in three spatial dimensions we have in total $3+2+1=6$ phase space variables i.e. 3 spatial ($x\in G$), 2 angular ($\omega\in S$) and 1 energy dimensions ($E\in I$).
Hence the numerical dimension of the problem is very large
in the sense that the total number of grid points needed in any discretization scheme
grows fast 
with the number of grid points used for discretizing each individual dimension.
The numerical computations can be carried out by utilizing e.g. diffusion approximation, spherical harmonics, finite element methods (FEM) or moment methods. 
Recently several contributions for radiotherapy purposes applying transport equations have been reported.
For forward problems (dose calculation) we refer  e.g. to \cite{larsen}, \cite{duclous}, \cite{hensel}, \cite{olbrant}, \cite{schneider}, \cite{pichard}, \cite{vassiliev} and both  forward and inverse problems (optimal control) to \cite{bomanthesis},
\cite{barnard},  \cite{frank08},   \cite{frank10},   \cite{herty11},  \cite{jorres}, \cite{tervo18}. In part of the  contributions  also numerical simulations have been exposed.

Since the object function (\ref{io}) contains non-convex terms, {\it global optimization} (\cite{pinter}) for  Lipschitz continuous functions in convex domains is needed. Moreover, the applied optimization method should be reasonably fast.
The initialization (that is, the determination of an {\it initial solution} for global optimization scheme) is necessary since the determination of a carefully chosen initial point for a large dimensional global optimization 
scheme is very essential for achieving time saving and satisfactory results (\cite{pinter14}).
We prove in section \ref{cis}  that for a certain (related) convex object function,
the optimal control exists, and we give formulas  for it in a variational form.
We suggest that this solution can be used as an initial solution for the actual global optimization scheme.

\section{Notations and Some Preliminary Analysis}\label{prem}

\subsection{Basic notations}\label{b-f-spaces}

We assume that $G$ is an open bounded 
 set in $\R^3$ such that $\ol{G}$ is a $C^1$-manifold with boundary \cite{lee}.
In particular, it follows from this definition that $G$ lies on one side of its boundary.
The unit outward (with respect to $G$) pointing normal on $\partial G$ is denoted by $\nu$ 
and the surface measure (induced by the Lebesgue measure $dx$) on $\partial G$ is denoted as $d\sigma$.
We let $S=S_2$ be the unit sphere in $\R^3$ equipped with the standard surface measure $d\omega$. 
Furthermore, let $I=[E_0,E_{\rm m}]$  where
$0\leq E_0<E_{\rm m}<\infty$. 
We shall denote by $I^\circ$ the interior  of $I$ and $I$ is equipped  with the Lebesgue measure $dE$.
All functions considered in this paper are real-valued, and all linear Hilbert spaces are real.

For $(x,\omega)\in G\times S$ the \emph{escape time} $t(x,\omega)=t_-(x,\omega)$ is defined as 
$ t(x,\omega):=\inf\{s>0\ |\ x-s\omega\not\in G\}  =\sup\{T>0\ |\ x-s\omega\in G\ {\rm for\ all}\ 0<s<T\}.$
We define
\[
\Gamma':=(\partial G)\times S,\quad
\Gamma:=\Gamma'\times I,
\]
and their subsets
\begin{alignat*}{3}
\Gamma_0' :={}&\{(y,\omega)\in \Gamma'\ |\ \omega\cdot\nu(y)=0\},
\quad &&\Gamma_0:=\Gamma'_0\times I, \\
\Gamma_{-}':={}&\{(y,\omega)\in \Gamma'\ |\ \omega\cdot\nu(y)<0\}, && \Gamma_{-}:=\Gamma_{-}'\times I,\\
\Gamma_{+}':={}&\{(y,\omega)\in \Gamma'\ |\ \omega\cdot\nu(y)>0\}, &&\Gamma_{+}:=\Gamma_{+}'\times I.
\end{alignat*}
Note that $\Gamma=\Gamma_0\cup \Gamma_-\cup \Gamma_+$.

In the sequel we denote for $k\in\N_0$, 
$
C^k(\ol G\times S\times I):=\{\psi\in C^k(G\times S\times I^\circ)\ |\ \psi=f_{|G\times S\times I^\circ},\ f\in C_0^k(\R^3\times S\times\R)\},
$
where for a $C^k$-manifold $M$ without boundary, the set $C_0^k(M)$ denotes
 the set of all $C^k$-functions on $M$ with compact support.
Define the (Sobolev) space $W^2(G\times S\times I)$  by
\[
W^2(G\times S\times I)
:=\{\psi\in L^2(G\times S\times I)\ |\  \omega\cdot\nabla_x \psi\in L^2(G\times S\times I) \}.
\]
The  space $W^2(G\times S\times I)$ is a Hilbert space  equipped  with the standard inner product
\[
\la {\psi},v\ra_{W^2(G\times S\times I)}:=\la {\psi},v\ra_{L^2(G\times S\times I)}+
\la\omega\cdot\nabla_x\psi,\omega\cdot\nabla_x v\ra_{L^2(G\times S\times I)}.
\]
Moreover, the space $C^1(\ol G\times S\times I)$ is a dense subspace of  $W^2(G\times S\times I)$ (\cite{friedrich44}).

Let $T^2(\Gamma)$ be the weighted Lebesgue space $L^2(\Gamma,|\omega\cdot\nu|d\sigma d\omega dE)$.
The trace $\gamma(\psi):=\psi_{|\Gamma}$ 
is well-defined for $\psi\in W^2(G\times S\times I)$ and $\gamma(\psi)\in L^2_{\rm loc}(\Gamma_{\pm },|\omega\cdot\nu|d\sigma d\omega dE)$ (\cite{cessenat84}, \cite[Section 2.2]{tervo18-up}). The space
\[
\widetilde{W}^2(G\times S\times I):=\{\psi\in W^2(G\times S\times I)\ |\ \gamma(\psi)\in T^2(\Gamma) \}
\]
is a Hilbert space  equipped with the inner product
\[
\la {\psi},v\ra_{\widetilde{W}^2(G\times S\times I)}&:=\la {\psi},v\ra_{W^2(G\times S\times I)}+ \la {\gamma(\psi)},\gamma(v)\ra_{T^2(\Gamma)};
\;\\
& 
\la h_1,h_2\ra_{T^2(\Gamma)}:=\int_\Gamma h_1(x,\omega,E)h_2(x,\omega,E)|\omega\cdot\nu| d\sigma d\omega dE.
\]

We need the following additional Hilbert spaces. 
Let $H_1$ be the completion of $C^1(\ol G\times S\times I)$ with respect to the inner product
\be
\la \psi,v\ra_{H_1}:=\la \psi,v\ra_{L^2(G\times S\times I)}+
\la\gamma(\psi),\gamma(v)\ra_{T^2(\Gamma)}
\ee
The elements of $H_1$ are  of the form $\tilde\psi=(\psi,q)\in 
L^2(G\times S\times I)\times T^2(\Gamma)$. Actually, they are exactly 
elements of the closure of the graph of trace operator $\gamma:C^1(\ol G\times S\times I)\to 
C^1(\partial G\times S\times I)$ in $L^2(G\times S\times I)\times T^2(\Gamma)$. 
The inner product in $H_1$ is
\begin{align}\label{spaceH}
\la\tilde\psi,\tilde \psi'\ra_{H_1}=\la\psi, \psi'\ra_{L^2(G\times S\times I)}
+\la q,q'\ra_{T^2(\Gamma)},
\end{align}
for $\tilde\psi=(\psi,q)$, $\tilde\psi'=(\psi',q')\in H_1$.
Finally, let
$H_2$ be the completion of $C^1(\ol G\times S\times I)$ with respect to the inner product
\be\label{inph2}
\la \psi,v\ra_{H_2}
:=\la \psi,v\ra_{H_1}+
\la\omega\cdot\nabla_x \psi,\omega\cdot\nabla_x v\ra_{L^2(G\times S\times I)}.
\ee
Obviously $H_2\subset H_1$
and for $\psi\in H_2$ we have $q=\psi_{|\Gamma}$. The inner product in $H_2$ is given by (\ref{inph2}).

We recall the   Green's formula \cite[p. 225]{dautraylionsv6} for every $\psi,\ v\in \widetilde{W}^2( G\times S\times I)$
\be\label{green}
\int_{G\times S\times I}(\omega\cdot \nabla_x \psi)v\, dxd\omega dE
+\int_{G\times S\times I}\psi(\omega\cdot \nabla_x v)\, dxd\omega dE=
\int_{\partial G\times S\times I}(\omega\cdot \nu)  \psi v d\sigma d\omega dE.
\ee

\subsection{Hyper-singular transport operator for M\o ller scattering}\label{h-sing-m}

The M\o ller collision operator has been considered e.g. in  \cite{tervo18-up,tervo18}. It  can be written explicitly for $E\geq E_0$ as
\bea\label{mol-a}
&
(K\psi)(x,\omega,E)=
{\rm p.f.}\int_E^{E_m}\hat\sigma_{2}(x,E',E){1\over{(E'-E)^2}}
\int_{0}^{2\pi}\psi(x,\gamma(E',E,\omega)(s),E')dsdE'
\nonumber\\
&
-
{\rm p.f.}\int_E^{E_m}\hat\sigma_{1}(x,E',E){1\over{E'-E}}
\int_{0}^{2\pi}\psi(x,\gamma(E',E,\omega)(s),E')dsdE'
+(K_0\psi)(x,\omega,E).
\eea
Here $K_0$ is given by 
\be\label{psio-0b}
(K_0\psi)(x,\omega,E):=
\int_{I'}\int_0^{2\pi}\chi(E',E)\hat\sigma_0(x,E',E)\psi(x,\gamma(E',E,\omega)(s),E') ds dE'
\ee
where $\chi(E',E)=\chi_{\R_+}(E'-E)$ (here $\chi_{\R_+}$ is the characteristic function $\R_+$).
Moreover,
$\gamma=\gamma(E',E,\omega):[0,2\pi]\to S$
is a parametrization of the curve
\[
\Gamma(E',E,\omega)=\{\omega'\in S\ |\ \omega'\cdot\omega-\mu(E',E)=0\},\
\mu(E',E):=\sqrt{{{E(E'+2)}\over{E'(E+2)}}}.
\]
We choose
\[
\gamma(E',E,\omega)(s)=R(\omega)\big(\sqrt{1-\mu^2}\cos(s),\sqrt{1-\mu^2}\sin(s),\mu\big),\quad s\in [0,2\pi],
\]
where $R(\omega)$ is a rotation (unitary) matrix which maps the vector $e_3=(0,0,1)$ to $\omega$.

Denote as in \cite{tervo18-up}, \cite{tervo18}
\be\label{psio-0a}
(K_j\psi)(x,\omega,E):=
{\rm p.f.}\int_E^{E_m}\hat\sigma_{j}(x,E',E){1\over{(E'-E)^j}}
\int_{0}^{2\pi}\psi(x,\gamma(E',E,\omega)(s),E')dsdE',\ j=1,2
\ee
and
\bea\label{Kj}
(\ol {\s K}_{j}\psi)(x,\omega,E',E)
:=\hat{\sigma}_{j}(x,E',E)\int_{0}^{2\pi}\psi(x,\gamma(E',E,\omega)(s),E')ds,
\eea
Then for $j=1,2$
\[
(K_j\psi)(x,\omega,E)=
{\s H}_j\big((\ol{\s K}_j\psi)(x,\omega,\cdot,E)\big)(E) 
\]
where ${\s H}_j$ are the Hadamard finite part singular integral operators (\cite{tervo18-up}, section 3.2).
The corresponding transport operator is
\bea\label{exth-1}
&
T\psi=-
{\s H}_2\big((\ol{\s K}_2\psi)(x,\omega,\cdot,E)\big)(E) 
+
{\s H}_1\big((\ol{\s K}_1\psi)(x,\omega,\cdot,E)\big)(E)\nonumber\\
&
+\omega\cdot\nabla_x\psi+\Sigma\psi
-K_r\psi.
\eea
A refined pseudo-differential like form of the operator (\ref{exth-1}) is derived in \cite{tervoarxiv-18}, Theorem 4.17 (together with Theorem 4.9, Lemma 4.10, Remark 4.11 therein).

\subsection{CSDA-Boltzmann type approximation}\label{tr-app}

We recall from \cite{tervo19} the truncation of the hyper-singular kernels and approximations leading to relevant approximation of the exact transport operator (\ref{exth-1}).  
The operator $K_0$ given by (\ref{psio-0b}) is an ordinary partial Schur integral operator and so it suffices only to truncate the first and second terms $(K_j\psi)(x,\omega,E)$ of equation (\ref{mol-a}).
In the following we assume that the "cut-off energy for primary particles" is $E'=\kappa E$ where $\kappa>1$.
Then, we decompose the integration in (\ref{psio-0a}) as follows for $j=1,2:$
\bea\label{psio-0c}
&
(K_j\psi)(x,\omega,E)=
\nonumber\\
&
=
{\rm p.f.}\int_E^{\kappa E}\hat\sigma_{j}(x,E',E){1\over{(E'-E)^j}}
\int_{0}^{2\pi}\psi(x,\gamma(E',E,\omega)(s),E')dsdE'\nonumber\\
&
+
\int_{\kappa E}^{E_m}\hat\sigma_{j}(x,E',E){1\over{(E'-E)^j}}
\int_{0}^{2\pi}\psi(x,\gamma(E',E,\omega)(s),E')dsdE'  \\[2mm]
& =: (K_{j,1,\kappa}\psi)(x,\omega,E) + (K_{j,0,\kappa}\psi)(x,\omega,E) \label{de}
\eea
where we noticed that the last integral $(K_{j,0,\kappa}\psi)$ does {\bf  not} involve a singular kernel.

Using these notations, we see by (\ref{mol-a}) that
\be\label{de-T}
T\psi=-K_{2,1,\kappa}\psi+K_{1,1,\kappa}\psi+\omega\cdot\nabla_x\psi+\Sigma\psi-(K_{2,0,\kappa}\psi-K_{1,0,\kappa}\psi+K_0\psi).
\ee
Since
\bea\label{psio-0d}
&
(K_{j,0,\kappa}\psi)(x,\omega,E)\nonumber\\
&
=
\int_{I'}\chi_{\R_+}(E'-\kappa E)\hat\sigma_{j}(x,E',E){1\over{(E'-E)^j}}
\int_{0}^{2\pi}\psi(x,\gamma(E',E,\omega)(s),E')dsdE'
\eea
we  find   that the operators $K_{j,0,\kappa}$ are partial Schur integral operators and hence they are bounded operators $L^2(G\times S\times I)\to L^2(G\times S\times I)$ (see section \ref{res-coll} below). 
That is why it suffices to consider approximations only for the operators 
$K_{j,1,\kappa}$.

Let for $j=1,2$
\[
f_j(x,\omega,E',E):=(\ol {\s K}_{j}\psi)(x,\omega,E',E)=
\hat\sigma_{j}(x,E',E)
\int_{0}^{2\pi}\psi(x,\gamma(E',E,\omega)(s),E')ds.
\]
For $K_{1,1,\kappa}$ we apply the approximation
\be\label{trunc-2-f1}
f_1(x,\omega,E',E)\approx  f_1(x,\omega,E,E)
\ee
\emph{on the interval $[E,\kappa E]$} which leads  to the approximation
\be\label{trunc-1}
(K_{1,1,\kappa}\psi)(x,\omega,E)\approx (\widetilde K_{1,1,\kappa}\psi)(x,\omega,E):=
2\pi\ \ln(\kappa E-E)\hat\sigma_{1}(x,E,E)\psi(x,\omega,E)
\ee
where we recalled that $\gamma(E,E,\omega)(s)=\omega$ and that (\cite{tervo18-up}, section 3.1)
\[
{\rm p.f.}\int_E^{\kappa E}{1\over{E'-E}}dE'= \ln(\kappa E-E).
\]

For $K_{2,1,\kappa}$ we apply the approximation
\be\label{trunc-2}
f_2(x,\omega,E',E)\approx  f_2(x,\omega,E,E)+{\p {f_2}{E'}}(x,\omega,E,E)(E'-E)
\ee
\emph{on the interval $[E,\kappa E]$}.  Then it can be shown that 
(see \cite{tervoarxiv-18} Theorem 4.9, Lemma 4.10 and Remark 4.11) 
\[
&
{\p {f_2}{E'}}(x,\omega,E,E)
=
2\pi\ {\p {\hat\sigma_{2}}{E'}}(x,E,E)\psi(x,\omega,E)\\
&
+\hat\sigma_{2}(x,E,E)(A(x,\omega,E,\partial_\omega)\psi)(x,\omega,E)
+2\pi\ \hat\sigma_{2}(x,E,E){\p \psi{E}}(x,\omega,E).
\]
Here
\be\label{aa}
A(\omega,E,\partial_\omega)\psi=
-\pi\ (\partial_{E'}\mu)(E,E)\Delta_S\psi
\ee
where $\Delta_S$ is the Laplace-Beltrami operator on $S$.
Hence we obtain the approximation
\bea\label{trunc-4}
&
(K_{2,1,\kappa}\psi)(x,\omega,E)
\approx
(\widetilde K_{2,1,\kappa}\psi)(x,\omega,E)\nonumber\\
&
:=
-2\pi\ {1\over{\kappa E-E}}\hat\sigma_{2}(x,E,E)\psi(x,\omega,E) \nonumber\\
&
+\ln(\kappa E-E)\hat\sigma_{2}(x,E,E)(A(x,\omega,E,\partial_\omega)\psi)(x,\omega,E)
+2\pi\ \hat\sigma_{2}(x,E,E)\ln(\kappa E-E){\p \psi{E}}(x,\omega,E)\nonumber\\
&
+2\pi\ \ln(\kappa E-E){\p {\hat\sigma_{2}}{E'}}(x,E,E)\psi(x,\omega,E)
\eea
where we again recalled that $\gamma(E,E,\omega)(s)=\omega$ 
and that 
\[
{\rm p.f.}\int_E^{\kappa E}{1\over{E'-E}}dE'=\ln(\kappa E-E),\
{\rm p.f.}\int_E^{\kappa E}{1\over{(E'-E)^2}}dE'=-{1\over{\kappa E-E}}.
\]

We introduce the following definitions,
\begin{align*}
S_\kappa(x,E)& := 2\pi\ \hat\sigma_{2}(x,E,E)\ln(\kappa E-E),
\\
\Sigma_\kappa(x,E)& :=\Sigma(x,E)+2\pi\ \hat\sigma_{2}(x,E,E){1\over{\kappa E-E}}\\
& -2\pi\ \ln(\kappa E-E){\p {\hat\sigma_{2}}{E'}}(x,E,E) +2\pi\ \ln(\kappa E-E)\hat\sigma_{1}(x,E,E),
\\
(Q_\kappa(x,\omega,E,\partial_\omega)\psi)(x,\omega,E)& :=
\ln(\kappa E-E)\hat\sigma_{2}(x,E,E)(A(x,\omega,E,\partial_\omega)\psi)(x,\omega,E),
\\
K_{r,\kappa}:=& K_0+K_{2,0,\kappa}-K_{1,0,\kappa}.
\end{align*}
Then,  the truncated and approximated transport operator, say $T_\kappa$, is given by (recall (\ref{de-T}), (\ref{trunc-1}), (\ref{trunc-4}))
\bea\label{trunc-6}
&
T_\kappa\psi:=-\widetilde K_{2,1,\kappa}\psi+\widetilde K_{1,1,\kappa}\psi
+\omega\cdot\nabla_x\psi+\Sigma\psi-K_{r,\kappa}\psi
\nonumber\\
&
=
-S_\kappa(x,E){\p \psi{E}}
-Q_\kappa(x,\omega,E,\partial_\omega)\psi
+\omega\cdot\nabla_x\psi+\Sigma_\kappa\psi-K_{r,\kappa}\psi.
\eea
Combining the above formulas we see that 
\[
Q_\kappa(x,\omega,E,\partial_\omega)\psi=
-\pi\ln(\kappa E-E)\hat\sigma_2(x,E,E)(\partial_{E'}\mu)(E,E)\Delta_S\psi.
\]

The restricted collision operator $K_{r,\kappa}$ is
\bea\label{K-r}
& 
(K_{r,\kappa}\psi)(x,\omega,E)=((K_0+K_{2,0,\kappa}-K_{1,0,\kappa})\psi)(x,\omega,E)\nonumber\\
&
=
\int_{I'}\chi_{\R_+}(E'-E)\hat\sigma_{0}(x,E',E)
\int_{0}^{2\pi}\psi(x,\gamma(E',E,\omega)(s),E')dsdE'\nonumber\\
&
-
\int_{I'}\chi_{\R_+}(E'-\kappa E)\hat\sigma_{1}(x,E',E){1\over{E'-E}}
\int_{0}^{2\pi}\psi(x,\gamma(E',E,\omega)(s),E')dsdE'\nonumber\\
&
+\int_{I'}\chi_{\R_+}(E'-\kappa E)\hat\sigma_{2}(x,E',E){1\over{(E'-E)^2}}
\int_{0}^{2\pi}\psi(x,\gamma(E',E,\omega)(s),E')dsdE'\nonumber\\
&
=
\int_{I'}\hat\sigma_{r,\kappa}(x,E',E)
\int_{0}^{2\pi}\psi(x,\gamma(E',E,\omega)(s),E')dsdE'
\eea
where
\[
&
\hat\sigma_{r,\kappa}(x,E',E):=
\chi_{\R_+}(E'-E)\hat\sigma_{0}(x,E',E)
-\chi_{\R_+}(E'-\kappa E)\hat\sigma_{1}(x,E',E){1\over{E'-E}}\\
&
+
\chi_{\R_+}(E'-\kappa E)\hat\sigma_{2}(x,E',E){1\over{(E'-E)^2}}.
\]

In \cite{tervo19} we verified that (in a certain sense) the approximative operator (\ref{trunc-6}) convergences to  that exact operator (\ref{exth-1})
as $\kappa$ approaches to $1$.

\section{Existence of solutions for the single transport problem }\label{trunc-exists}

\subsection{Assumptions for the restricted collision operator}\label{res-coll}

In  single charged particle transport problem we take the below assumptions
for
the restricted collision operator $K_r$. A realistic transport model comprises a coupled system of transport equations which will be considered in the next section \ref{cosyst}. Therein we also introduce a somewhat more general  restricted collision operator.

We assume that
\be\label{coll-3}
(K_r\psi)(x,\omega,E)
=
\int_{I'}\int_{0}^{2\pi}
\hat\sigma_r(x,E',E)
\psi(x,\gamma(E',E,\omega)(s),E')ds dE'
\ee
where
$\hat\sigma_r:G\times I'\times I\to\R$ is a non-negative measurable function such that
\bea\label{ass-8b}
&
\int_{I'}\hat{\sigma}_r(x,E',E)dE'\leq M_1,  \ {\rm a.e.}\ (x,E)\in G\times I, \nonumber\\
&
\int_{I'}\hat{\sigma}_r(x,E,E')dE'\leq M_2,  \ {\rm a.e.}\ (x,E)\in G\times I
\eea
Moreover, we assume that
\be\label{ass8-aa}
\Sigma(x,\omega,E)
-2\pi\int_{I'}\hat{\sigma}_r(x,E,E')dE'
\geq c,
\ee
and
\be\label{ass9-a}
\Sigma(x,\omega,E)
-2\pi\int_{I'}\hat{\sigma}_r(x,E',E)dE'
\geq c.
\ee
for a.e. $(x,\omega,E)\in G\times S\times I$ where $c>0$.

We recall from \cite{tervo18}

\begin{theorem}\label{SK-dissip}
Suppose that 
$\Sigma\in L^\infty(G\times S\times I),\ \Sigma\geq 0
$
and that
the assumptions   (\ref{ass-8b}), (\ref{ass8-aa}) and (\ref{ass9-a}) are valid. 
Then $K_r$  is a bounded operator 
$L^2(G\times S\times I)\to L^2(G\times S\times I)$ and 
\be\label{coer}
\la (\Sigma-K_r)\psi,\psi\ra_{L^2(G\times S\times I)}\geq c\n{\psi}^2_{L^2(G\times S\times I)}\ {\rm for\ all}\ \psi\in L^2(G\times S\times I).
\ee
 
\end{theorem}

\begin{proof}  \end{proof}
 
\begin{remark}\label{srk}

Note that the integrals appearing in the  conditions  (\ref{ass-8b}), 
(\ref{ass8-aa}), (\ref{ass9-a})  
for the above approximative operator $K_{r,\kappa}$ (that is, for (\ref{K-r})) become 
\[
&
\int_{I'}\hat\sigma_{r,\kappa}(x,E',E) dE'
\\
&
=
\int_{E}^{E_m} \hat\sigma_{0}(x,E',E)dE'
-\int_{\kappa E}^{E_m}\hat\sigma_{1}(x,E',E){1\over{E'-E}}dE'\\
&
+
\int_{\kappa E}^{E_m}\hat\sigma_{2}(x,E',E){1\over{(E'-E)^2}} 
dE'
\]
and
\[
&
\int_{I'}\hat\sigma_{r,\kappa}(x,E,E') dE'
\\
&
=
\int_{E_0}^{E} \hat\sigma_{0}(x,E,E')dE'
-\int_{E_0}^{E/\kappa}\hat\sigma_{1}(x,E,E'){1\over{E-E'}}dE'\\
&
+
\int_{E_0}^{E/\kappa}\hat\sigma_{2}(x,E,E'){1\over{(E-E')^2}} 
dE'.
\]
 
\end{remark}

\subsection{Variational formulations for a single transport problem}\label{var-form}

The transport problem contains  solving an initial inflow boundary value problem of the form
\be\label{trunc-ex-1}
T\psi=f,\ \psi_{|\Gamma_-}=g,\ \psi(.,.,E_m)=0
\ee
where $T$ is an  approximation of the exact transport operator and where $f\in L^2(G\times S\times I),\ g\in T^2(\Gamma_-)$. 
The appropriate approximative operator $T=T_\kappa$ is given, by the discussion in the previous section,  by equation \eqref{trunc-6}. 
In the following we assume that 
$T$ is of the form, 
\be\label{trunc-ex-2}
T\psi=
a(x,E){\p \psi{E}}
+b(x,E)\Delta_S\psi
+\omega\cdot\nabla_x\psi+\Sigma\psi-K_{r}\psi.
\ee
Note that the operator $T_\kappa$ in (\ref{trunc-6}) is of the form (\ref{trunc-ex-2}).  

\begin{remark}\label{re-ge} 
We could replace $b(x,E)\Delta_S\psi$ with a more general operator
\be\label{b-type}
b(x,\omega,E,\partial_\omega)\psi
= b(x,E)\Delta_S\psi
+d(x,\omega,E)\cdot\nabla_S\psi.
\ee
where
$d=(d_1,d_2)\sim d_1{\partial\over{\partial\omega_1}}+ d_2{\partial\over{\partial\omega_2}}$ and $d\cdot\nabla_S$ is the Riemannian inner product on $T(S)$. After minor changes the below considerations are valid for this more general setting.
\end{remark}

At first
we verify formally the corresponding variational equation.
Assume that $\psi$ is a solution of (\ref{trunc-ex-1}) in the classical sense.
Let $\psi,\ v\in C^1(\ol G\times I,C^2(S))$.
By the  properties of Laplace-Beltrami operator we have 
\[
\int_S(\Delta_S\psi)(x,\omega,E)v(x,\omega,E) d\omega=
-\int_S\la (\nabla_S\psi)(x,\omega,E),(\nabla_Sv)(x,\omega,E)\ra d\omega 
\]
where $\la \nabla_S\psi,\nabla_Sv\ra$ is the chosen Riemannian inner product on $S$. Hence, 
\bea\label{add-a} 
&
\la b(x,E)\Delta_S\psi,v\ra_{L^2(G\times S\times I)}
=
- \la \nabla_S\psi,b(x,E)
\nabla_Sv\ra_{L^2(G\times S\times I)}\nonumber\\
&
:=-\int_G\int_I\int_S\la (\nabla_S\psi)(x,\omega,E),b(x,E)(\nabla_Sv)(x,\omega,E)\ra d\omega dE dx.
\eea
 
Using  integration by parts we have
\bea\label{trunc-9}
&
\la a{\p \psi{E}},v\ra_{L^2(G\times S\times I)}
=\int_G\int_Sa(x,E_m)\psi(x,\omega,E_m)v(x,\omega,E_m)d\omega dx\nonumber\\
&
-\int_G\int_Sa(x,E_0)\psi(x,\omega,E_0)v(x,\omega,E_0)d\omega dx
-\la \psi,{\p {(av)}{E}}\ra_{L^2(G\times S\times I)} \nonumber\\
&
=
-
\la \psi,a {\p v{E}}\ra_{L^2(G\times S\times I)} 
-\la \psi,{\p {a}{E}}v\ra_{L^2(G\times S\times I)}
+\la a(.,E_m)\psi(.,.,E_m),v(.,.,E_m)\ra_{L^2(G\times S)}\nonumber\\
&
-\la a(.,E_0)\psi(.,.,E_0),v(.,.,E_0)\ra_{L^2(G\times S)}
.
\eea
Moreover, by the  Green's formula (\ref{green}) we obtain 
\bea\label{trunc-10}
&
\la\omega\cdot\nabla_x\psi,v\ra_{L^2(G\times S\times I)}=
-\la\psi,\omega\cdot\nabla_x v\ra_{L^2(G\times S\times I)}
+\int_{\partial G\times S\times I}(\omega\cdot\nu)\psi v d\sigma d\omega dE  
\nonumber\\
&
=
-\la\psi,\omega\cdot\nabla_x v\ra_{L^2(G\times S\times I)}
+\int_{\Gamma_+}(\omega\cdot\nu)_+\psi v d\sigma d\omega dE
-\int_{\Gamma_-}(\omega\cdot\nu)_-\psi v d\sigma d\omega dE  
\eea
where $(\omega\cdot\nu)_{\pm}$ is the positive/negative part of $\omega\cdot\nu$.
Finally,
\be\label{trunc-12}
\la K_{r}\psi,v\ra_{L^2(G\times S\times I)}=
\la \psi, K_{r}^*v\ra_{L^2(G\times S\times I)},\ 
\la \Sigma\psi,v\ra_{L^2(G\times S\times I)}=
\la \psi, \Sigma^*v\ra_{L^2(G\times S\times I)}
.
\ee

As a conclusion we see that if $\psi$ is a classical solution of problem (\ref{trunc-ex-1}) then the following weak formulation is fulfilled
\bea\label{trunc-11}
&
B(\psi,v):=
-\la \psi,a {\p v{E}}\ra_{L^2(G\times S\times I)} 
-\la \psi,{\p {a}{E}}v\ra_{L^2(G\times S\times I)} \nonumber\\
&
-\la a(.,E_0)\psi(.,.,E_0),v(.,.,E_0)\ra_{L^2(G\times S)}\nonumber\\
&
- \la \nabla_S\psi,b(x,E)
\nabla_Sv\ra_{L^2(G\times S\times I)}
\nonumber\\
&
-\la\psi,\omega\cdot\nabla_xv\ra_{L^2(G\times S\times I)}
+\la\gamma_+(\psi),\gamma_+(v)\ra_{T^2(\Gamma_+)}
\nonumber\\
&
+
\la \psi, \Sigma^* v -K_{r}^*v\ra_{L^2(G\times S\times I)}
=\la f,v\ra_{L^2(G\times S\times I)}+\la g,\gamma_-(v)\ra_{T^2(\Gamma_-)}.
\eea

Define in $C^1(\ol G\times I,C^2(S))$  inner products
\bea \label{inH0}
&
\la \psi,v\ra_{\mc H}:=\la \psi,v\ra_{L^2(G\times S\times I)}+
\la\gamma(\psi),\gamma(v)\ra_{T^2(\Gamma)}\nonumber\\
&
+\la\psi(.,.,E_0),v(.,.,E_0)\ra_{L^2(G\times S)}+\la\psi(.,.,E_m),v(.,.,E_m)\ra_{L^2(G\times S)}\nonumber\\
&
+
\la \psi,v\ra_{L^2(G\times I,H^1(S))}
\eea
and 
\be\label{inhatH0}
\la \psi,v\ra_{\widehat {\mc H}}
=\la \psi,v\ra_{H}
+
\la\omega\cdot\nabla_x \psi,\omega\cdot\nabla_x v\ra_{L^2(G\times S\times I)}
+
\la {\p \psi{E}},{\p v{E}}\ra_{L^2(G\times S\times I)}.
\ee
Let $\mc H$ and $\widehat{\mc H} $ be the completions of
$C^1(\ol G\times I, C^2(S))$ with respect to inner products $\la .,.\ra_{\mc H}$ and $\la .,.\ra_{\mc {\widehat H}}$, respectively.

We assume for the coefficients now:
\be\label{trunc-18}
a\in L^\infty(G,W^{\infty,1}(I)), b\in L^\infty(G\times I),
\ee
\be\label{trunc-14}
-{\p a{E}}(x,E)\geq q_1>0, \ {\rm a.e.},
\ee
\be\label{trunc-19}
-b(x,E)\geq q_2>0,  \ {\rm a.e.},
\ee
\be\label{trunc-20}
-a(x,E_0)\geq q_3>0,\ -a(x,E_m)\geq q_3>0, \  {\rm a.e.}.
\ee

In \cite{tervo19}
we showed that under these assumptions the bilinear form $B:C^1(\ol G\times I,C^2(S))\times C^1(\ol G\times I,C^2(S))\to\R$ obeys the following {\it boundedness and coercivity} results:

\begin{theorem}\label{csdath1}
Suppose that  the assumptions (\ref{ass-8b}), (\ref{ass8-aa}), (\ref{ass9-a}),
(\ref{trunc-18}), (\ref{trunc-14}), (\ref{trunc-19}), (\ref{trunc-20}) are valid.
Then there exists a constant $M>0$  such that 
\be\label{csda29}
|B(\psi,v)|\leq M\n{\psi}_{\mc H}\n{v}_{\widehat{\mc H}}\quad \forall\psi,\ v\in C^1(\ol G\times  I,C^2(S))
\ee
and 
\be\label{csda30}
B(v,v)\geq c'\n{v}_{\mc H}^2
\quad \forall v\in C^1(\ol G\times I,C^2(S))
\ee
where 
\bea\label{cprime}
c':=\min\{{{q_1}\over 2},{{q_3}\over 2},q_2, {1\over 2},c\}>0 .
\eea 
\end{theorem}

\begin{proof}  \end{proof}

Because $C^1(\ol G\times I,C^2(S))\times C^1(\ol G\times I,C^2(S))$ is dense in $\mc H\times \widehat{\mc H}$ and because \eqref{csda29} holds, the bilinear form $B(\cdot,\cdot): C^1(\ol G\times I,C^2(S))\times C^1(\ol G\times I,C^2(S))\to\R$ has a unique extension $\widetilde B(\cdot,\cdot):\mc H\times \widehat{\mc H}\to\R$
which satisfies 
\be\label{csda36}
|\widetilde B(\psi,v)|\leq M\n{\psi}_{\mc H}\n{v}_{\widehat {\mc H}}\quad \forall\ \psi\in \mc H,\ v\in \widehat{\mc H}
\ee
and 
\be\label{csda37}
\widetilde B(v,v)\geq c'\n{v}_{\mc H}^2\quad \forall v\in \widehat {\mc H}.
\ee
We see that actually
\bea\label{coex}
&
\widetilde B(\psi,v)=
-\la \psi,a {\p v{E}}\ra_{L^2(G\times S\times I)} 
-\la \psi,{\p {a}{E}}v\ra_{L^2(G\times S\times I)} \nonumber\\
&
-
\la a(.,E_0) p_0,v(.,.,E_0)\ra_{L^2(G\times S)}\nonumber\\
& 
-\la \nabla_S\psi,b(x,E)
\nabla_Sv\ra_{L^2(G\times S\times I)}
\nonumber\\
&
-\la\psi,\omega\cdot\nabla_xv\ra_{L^2(G\times S\times I)}
+\la q_{|\Gamma_+},\gamma_+(v)\ra_{T^2(\Gamma_+)}
\nonumber\\
&
+
\la \psi, \Sigma^* v -K_{r}^*v\ra_{L^2(G\times S\times I)}.
\eea
where $q\in T^2(\Gamma)$ and $p_0\in L^2(G\times S)$ are explained in \cite{tervo18-up}, p. 14 (here we have $E_0$ instead of $0$ therein).

In addition, since for $\psi\in C^1(\ol G\times I,C^2(S))$ 
\bea\label{csda39}
|F(\psi)|
\leq{}& |\la {f},\psi\ra_{L^2(G\times S\times I)}|+|\la { g},\gamma_-(\psi)\ra_{T^2(\Gamma_-)}|
\nonumber\\
\leq{}&
\n{{ f}}_{L^2(G\times S\times I)}\n{\psi}_{L^2(G\times S\times I)}+\n{{g}}_{T^2(\Gamma_-)}\n{\gamma(\psi)}_{T^2(\Gamma_-)},
\eea
the linear form $F:C^1(\ol{G}\times I,C^2(S))\to\R$ has a unique bounded extension,
which we still denote by $F$,
\begin{align}\label{Fex}
F:\mc H\to\R;\quad
F(\psi)=\la { f},\psi\ra_{L^2(G\times S\times I)}+\la { g}, q_{|\Gamma_-}\ra_{T^2(\Gamma_-)},
\end{align}
Note also that the embedding $\widehat {\mc H}\subset \mc H$ is continuous.

\begin{example}\label{kappa=2}

Consider the approximation (\ref{trunc-6}) for $\kappa=2$ (which is the usual choice in practise). We assume that $E_0>1$ and so $\ln(\kappa E-E)=\ln(E)\geq \ln(E_0)>0$.
In certain applications (e.g. in dose calculation of radiation therapy for primary electrons) the expressions $\hat\sigma_2(x,E',E)$ appearing in the above differential cross sections are of the form
\bea\label{hsigma}
&
\hat\sigma_2(x,E',E)=\sigma_{0}(x){{(E'+1)^2}\over{E'(E'+2)}},\\
&
\hat\sigma_1(x,E',E)=\sigma_{0}(x){{(E'+1)^2}\over{E'(E'+2)}}{{2E'+1}\over{(E'+1)^2E}} 
\\
&
\hat\sigma_0(x,E',E)=\sigma_{0}(x){{(E'+1)^2}\over{E'(E'+2)}}\Big({1\over{E^2}}
+{1\over{(E'+1)^2}}\Big).
\eea
We assume that $\sigma_0\in L^\infty(G)$ and
\[
\sigma_0(x)\geq q_0>0\ {\rm a.e.}.
\]

We have
\[
{\p {\hat\sigma_2}{E'}}(x,E',E)
=\sigma_0(x){{2(E'+1)E'(E'+2)-(E'+1)^2(2E'+2)}\over{E'^2(E'+2)^2}}
\]
and so
\bea\label{der-s2}
&
{\p {\hat\sigma_2}{E'}}(x,E,E)
=\sigma_0(x){{2(E+1)E(E+2)-(E+1)^2(2E+2)}\over{E^2(E+2)^2}}
\nonumber\\
&
=2\sigma_0(x){{(E+1)(E^2+2E-(E+1)^2)}\over{E^2(E+2)^2}}=-2\sigma_0(x){{E+1}\over{E^2(E+2)^2}}.
\eea
Note also that
\[
{\p {\hat\sigma_2}{E}}(x,E',E)=0.
\]

For
$
\mu(E',E)=\sqrt{{E(E'+2)}\over{E'(E+2)}}
$
we get
\[
&
2\mu(E',E)\partial_{E'}\mu(E',E)={{EE'(E+2)-(E+2)E(E'+2)}\over{E'^2(E+2)^2}}
\\
&
={{E^2E'+2EE'-E'E^2-2E^2-2EE'-4E}\over{E'^2(E+2)^2}}=-2{{E(E+2)}\over{E'^2(E+2)^2}}.
\]
Hence
\be\label{der-mu}
\partial_{E'}\mu(E,E)=-{1\over{E(E+2)}}.
\ee

We see that (for $\kappa=2$)
\[
-a(x,E)=-(-S_\kappa(x,E))
=2\pi\ \hat\sigma_{2}(x,E,E)\ln(E)
\geq q_3>0,
\]
and since $-\ln(E)\geq -E$ for $E\geq 1$ we have
\[
&
-{\p a{E}}(x,E)=-({\p {(-S_\kappa)}{E}})(x,E)\\
&
=2\pi\big({\p {\sigma_2}{E}}(x,E,E)+{\p {\sigma_2}{E'}}(x,E,E)\big)\ln(E)+2\pi\ \hat\sigma_{2}(x,E,E){1\over E}\\
&
=2\pi\sigma_0(x)\Big(-{{2(E+1)}\over{E^2(E+2)^2}}\ln(E)+{{(E+1)^2}\over{E(E+2)}}{1\over E}\Big)
\\
&
\geq
{{2\pi\sigma_0(x)}\over{E^2(E+2)^2}}((E+1)^2(E+2)-2(E+1)E)\\
&
={{2\pi\sigma_0(x)}\over{E^2(E+2)^2}}(E+1)((E+1)(E+2)-2E)
\geq q_1>0.
\]
Moreover,
\[
&
-b(x,E)=-\pi\ln(E)\hat\sigma_2(x,E,E)(\partial_{E'}\mu)(E,E)
\\
&
=\pi\hat\sigma_2(x,E,E){1\over{E(E+2)}}\ln(E)\geq q_2>0.
\]

Furthermore,
\[
&
\Sigma_\kappa(x,E) =\Sigma(x,E)+2\pi\ \hat\sigma_{2}(x,E,E){1\over{E}}\\
&
-2\pi\ \ln(E){\p {\hat\sigma_{2}}{E'}}(x,E,E) +2\pi\ \ln(E)\hat\sigma_{1}(x,E,E)
\]
and so $\Sigma_\kappa\in L^\infty(G\times S\times I)$ when $\Sigma\in L^\infty(G\times S\times I)$. The conditions (\ref{ass-8b}) can be easily seen to be valid. The conditions (\ref{ass8-aa}), (\ref{ass9-a}) mean that (recall Remark \ref{srk})
\bea\label{ex-k-1}
&
\Sigma_\kappa(x,E)-
\int_{E}^{E_m} \hat\sigma_{0}(x,E',E)dE'
+\int_{\kappa E}^{E_m}\hat\sigma_{1}(x,E',E){1\over{E'-E}}dE'\\
&
-
\int_{\kappa E}^{E_m}\hat\sigma_{2}(x,E',E){1\over{(E'-E)^2}} 
dE'\geq c>0
\eea
and
\bea\label{ex-k-2}
&
\Sigma_\kappa(x,E)-
\int_{E_0}^{E} \hat\sigma_{0}(x,E,E')dE'
+\int_{E_0}^{E/\kappa}\hat\sigma_{1}(x,E,E'){1\over{E-E'}}dE'\\
&
-
\int_{E_0}^{E/\kappa}\hat\sigma_{2}(x,E,E'){1\over{(E-E')^2}} 
dE'\geq c>0.
\eea
We omit further proceedings of these conditions.
\end{example}

\subsection{Existence of solutions of the  initial and boundary value  problem  for a single transport problem}\label{bo-in-pr}

Let
\[
P(x,\omega,E,D)\psi:= 
a(x,E){\p \psi{E}}
+b(x,E)\Delta_S\psi
+ \omega\cdot\nabla_x\psi.
\]
The space
\begin{multline}
{\s H}_P(G\times S\times I^\circ):=\{\psi\in L^2(G\times S\times I)\ | \\
P(x,\omega,E,D)\psi\in L^2(G\times S\times I)\ {\rm in\ the \ weak\ sense}\}
\end{multline}
is a Hilbert space when equipped with the inner product 
\[
\la \phi,v\ra_{{\s H}_P(G\times S\times I^\circ)}=\la \psi,v\ra_{L^2(G\times S\times I)}+\la P(x,\omega,E,D)\psi,P(x,\omega,E,D)v\ra_{L^2(G\times S\times I)}.
\]
With this notation, the equation \eqref{trunc-ex-1}  can be written as
\begin{align*}
P(x,\omega,E,D)\psi+\Sigma\psi - K_r\psi = f.
\end{align*}

For the first instance  we recall an existence of weak solutions \emph{without boundary and initial conditions} 

\begin{theorem}\label{csdath3}
Suppose that the assumptions  (\ref{ass-8b}), (\ref{ass8-aa}), (\ref{ass9-a}),
(\ref{trunc-18}), (\ref{trunc-14}), (\ref{trunc-19}), (\ref{trunc-20}) are valid.
Let ${ f}\in L^2(G\times S\times I)$ and $g\in T^2(\Gamma_-)$.
Then the  variational equation 
\be\label{csda40a}
\widetilde{B}(\psi,v)=F(v)\quad \forall v\in \widehat{\mc H}
\ee
has a solution $\psi\in \mc H$.
Furthermore, $\psi \in {\s H}_P(G\times S\times I^\circ)$ and it is a weak (distributional) solution of the equation (\ref{trunc-ex-1}) that is,
\be\label{trunc-ex-2-a}
T\psi:=a(x,E){\p \psi{E}}
+b(x,E)\Delta_S\psi
+\omega\cdot\nabla_x\psi+\Sigma\psi-K_{r}\psi=f.
\ee
 
\end{theorem}

\begin{proof}

See \cite{tervo19}, Theorem  6.8. 
 
\end{proof}

Let 
\[
Q(x,\omega,E,D)\psi=
a(x,E){\p \psi{E}}
+\omega\cdot\nabla_x\psi.
\]
Define the space ${\mathbb H}_Q(G\times S\times I^\circ)$ 
\[
{\mathbb H}_Q(G\times S\times I^\circ):=
\{\psi\in L^2(G\times I,H^1(S))|\ Q(x,\omega,E,D)\psi
\in L^2(G\times I,H^{-1}(S))\}
\]
equipped with the inner product
\[
\la\psi,v\ra_{{\mathbb H}_Q(G\times S\times I^\circ)}
=\la\psi,v\ra_{L^2(G\times I,H^1(S))}+\la Q(x,\omega,E,D)\psi, Q(x,\omega,E,D)v\ra_{L^2(G\times I,H^{-1}(S))}.
\]

Define the traces $\gamma_{\pm}(\psi):=\psi_{|\Gamma_{\pm}}$ and $\gamma_{\rm m}(\psi):=\psi(\cdot,\cdot,E_{\rm m}),\ \gamma_{0}(\psi):=\psi(\cdot,\cdot,E_0)$. 
The traces $\gamma_{\pm}:W(G\times S\times I)\to L_{\rm loc}^2(\Gamma_{\pm},|\omega\cdot\nu| d\sigma d\omega dE)$ can be shown to be continuous  (\cite{cessenat84} or \cite{tervo18-up}, Theorem 2.16). Moreover, we have

\begin{theorem}\label{exttr}
Suppose that $a\in  W^{\infty,1}(G\times I^\circ)$  such that 
\be\label{assfor-a}
|a(x,E)|\geq c_0>0\ {\rm a.e.\  in}\ \ol G\times I.
\ee
Then the trace operators
\begin{align*}
\gamma_{\pm}:{\mathbb H}_Q(G\times S\times I^\circ) & \to L_{\rm loc}^2(\Gamma_{\pm},|\omega\cdot\nu| d\sigma d\omega dE), \\
\gamma_{\rm m}:{\mathbb H}_Q(G\times S\times I^\circ) & \to L_{\rm loc}^2(G\times S), \\
\gamma_{0}:{\mathbb H}_Q(G\times S\times I^\circ) & \to L_{\rm loc}^2(G\times S),
\end{align*}
are well-defined and continuous.
\end{theorem}

\begin{proof}

See \cite{tervo19}, Theorem 3.10.

\end{proof}

Let $P^*(x,\omega,E,D)$ be the formal adjoint operator of $P(x,\omega,E,D)$ that is,
\[
P^*(x,\omega,E,D)v
=
-a {\p v{E}}
-{\p {a}{E}}v
+b(x,E)\Delta_Sv
-\omega\cdot\nabla_x v.
\]
For $U\in  L^2(G\times I,H^{-1}(S),\ v\in L^2(G\times I,H^{1}(S)$ we define
\[
\la U,v\ra:=\int_{G\times I}( U(x,\omega,.),v(x,\omega,.)) dx dE
\]
where $( U(x,\omega,.),v(x,\omega,.))$ is the canonical duality between $H^{-1}(S)$ and $H^1(S)$.
We have  the following extended Green's formula  

\begin{theorem}\label{ex-green}
Suppose that (\ref{trunc-18}), (\ref{assfor-a}) hold and that
$\psi\in {\s H}_P(G\times S\times I^\circ)\cap L^2(G\times I,H^{1}(S))$ and $v\in\widehat{\mc H}$ for which
$({\rm supp}(v))\cap \partial (G\times S\times I)$ is a compact subset of $\Gamma_-\cup \Gamma_+\cup (G\times S\times\{E_m\})\cup(G\times S\times\{E_0\})$. Then
\bea\label{green-ex}
&
\la P(x,\omega,E,D) \psi,v\ra
-\la P^*(x,\omega,E,D) v,\psi\ra 
=
\int_{\partial G\times S\times I}(\omega\cdot \nu) \psi v\ \ d\sigma d\omega dE\nonumber\\
&
+
\int_{G\times S}\big(a(\cdot,E_m)\psi(\cdot,\cdot,E_m)v(\cdot,\cdot,E_m)-a(\cdot,E_{0})\psi(\cdot,\cdot,E_{0})v(\cdot,\cdot,E_{0})\big)dx d\omega.
\eea

\end{theorem}

The proof follows by applying standard Green's formula for smooth functions and density arguments. The Theorem \ref{exttr} on the continuity of the trace operators  guarantees that the Green's formula  holds also in the limit. We omit details.

\begin{remark}\label{re-ex-green}

The Green formula has some additional generalizations. Especially, (\ref{green-ex}) holds for $\psi=v$ in the case when
$\psi\ \in {\s H}_P(G\times S\times I^\circ)\cap L^2(G\times I,H^{1}(S))$ 
such that $\gamma_{\pm}(\psi)\in T^2(\Gamma_{\pm})$ and $\gamma_{\rm m}(\psi),\ \gamma_0(\psi)\in L^2(G\times S)$ (cf. \cite{dautraylionsv6}, p. 225).  

\end{remark}

When (\ref{assfor-a}) holds
the weak solution $\psi$ of the
equation \eqref{trunc-ex-2-a} obtained in Theorem \ref{csdath3} can be shown to be a solution of the initial and inflow boundary value problem. 
We have

\begin{theorem}\label{cor-csdath3-a}
Suppose that the assumptions  (\ref{ass-8b}), (\ref{ass8-aa}), (\ref{ass9-a}),
(\ref{trunc-18}), (\ref{trunc-14}), (\ref{trunc-19}), (\ref{trunc-20}), (\ref{assfor-a}) are valid.
Let ${ f}\in L^2(G\times S\times I)$ and $g\in T^2(\Gamma_-)$.
Then the  transport problem 
\bea\label{trunc-ex-3-b}
&
a(x,E){\p \psi{E}}
+b(x,E)\Delta_S\psi
+\omega\cdot\nabla_x\psi+\Sigma\psi-K_{r}\psi=f\nonumber\\
&
\psi_{|\Gamma_-}=g,\ \psi(.,.,E_m)=0
\eea
has a unique solution $\psi\in \mc H\cap {\s H}_P(G\times S\times I^\circ)$.

In addition,   the solution $\psi$ obeys the apriori  estimate
\be\label{csda40aa}
\n{\psi}_{{\mc H}}\leq C
\big(\n{{ f}}_{L^2(G\times S\times I)}+\n{{g}}_{T^2(\Gamma_-)}\big).
\ee
\end{theorem}

\begin{proof}
The proof  follows by the techniques applied in the proof of Theorem 5.7, \cite{tervo18-up} (items (ii)-(iii) of the proof). We omit detailed treatments.
\end{proof}

\begin{remark}\label{g+}
Additionally, we recall  from \cite{tervo18-up}, Theorem 5.8 that under the assumptions of Theorem \ref{cor-csdath3-a}
$\psi_{|\Gamma_+}\in T^2(\Gamma_+)$ and  $\psi(.,.,E_0)\in L^2(G\times S)$.
\end{remark}


\section{Existence of Solutions for  Coupled Systems}\label{cosyst}

In reality e.g. the dose calculation model in radiation therapy is a coupled system of transport equations governing the simultaneous evolving of (at least) photons, electrons ans positrons. 
In this section, we consider the coupled transport problem.

\subsection{Coupled problem and existence of solutions}\label{coup-ex}

Let ${\bf f}=(f_1,f_2,f_3)\in L^2(G\times S\times I)^3$ and ${\bf g}=(g_1,g_2,g_3)\in T^2(\Gamma_-)^3$.
We deal with the following coupled system of partial integro-differential equations
for $\psi=(\psi_1,\psi_2,\psi_3)$ on $G\times S\times I$,
\begin{gather}
\omega\cdot\nabla_x\psi_1+\Sigma_1\psi_1-K_{r,1}\psi=f_1, \label{csda1a}\\
a_j(x,E){\p {\psi_j}{E}}
+b_j(x,E)\Delta_S\psi_j
+\omega\cdot\nabla_x\psi_j+\Sigma_{j}\psi_j-K_{r,j}\psi
=f_j,\quad j=2,3.\label{csda1b}
\end{gather}
In order to guarantee uniqueness of solutions, we moreover impose the inflow boundary conditions on $\Gamma_-$,
\be\label{csda2}
{\psi_j}_{|\Gamma_-}=g_j,\quad j=1,2,3,
\ee
and initial value  conditions on $G\times S$,
\be\label{csda3}
\psi_j(\cdot,\cdot,E_{\rm m})=0,\quad j=2,3,
\ee
where $E_{\rm m}$ is the cut-off energy.

We assume that 
\be\label{ass-S}
\Sigma_j\in L^\infty(G\times S\times I), \Sigma_j\geq 0\ {\rm a.e.}\ j=1,2,3.
\ee 
Define for $\psi=(\psi_1,\psi_2,\psi_3)\in L^2(G\times S\times I)^3$,
\be\label{sda1}
{\bf \Sigma}\psi=(\Sigma_1\psi_1,\Sigma_2\psi_2,\Sigma_3\psi_3)
\ee
and 
\be\label{sda2}
{\bf K}_r\psi=(K_{r,1}\psi,K_{r,2}\psi,K_{r,3}\psi).
\ee
We see that ${\bf \Sigma}:L^2(G\times S\times I)^3\to L^2(G\times S\times I)^3$ is a bounded linear operator.

Each $K_{r,j}$ may be a sum of different type of restricted collision operators 
\be\label{esols1} 
K_{r,j}=K_{r,j}^1+K_{r,j}^2+K_{r,j}^3
\ee
corresponding to different type of particle interactions; see section 5.4 of \cite{tervo18-up}.
Here $K_{r,j}^1$ is of the form 
\[
(K_{r,j}^1\psi)(x,\omega,E)=\int_{S'\times I'}\sigma_{r,j}^1(x,\omega',\omega,E',E)\psi(x,\omega',E')d\omega' dE',  
\]
where $\sigma_{r,j}^1:G\times S^2\times I^2\to\R$ is a non-negative measurable function such that 
\bea\label{ass-s-1}
&\int_{S'\times I'}\sigma_{r,j}^1(x,\omega',\omega,E',E)d\omega' dE'\leq M_1,\nonumber\\
&\int_{S'\times I'}\sigma_{r,j}^1(x,\omega,\omega',E,E')d\omega' dE'\leq M_2,
\eea
for a.e. $(x,\omega,E)\in G\times S\times I$.
The operator
$K_{r,j}^2$ is of the form 
\[
(K_{r,j}^2\psi)(x,\omega,E)=\int_{ S'}\sigma_{r,j}^2(x,\omega',\omega,E)\psi(x,\omega',E) d\omega',  
\]
where $\sigma_{r,j}^2:G\times S^2\times I\to\R$ is a non-negative  measurable function  such that
\bea\label{ass-s-2}
&\int_{S'}\sigma_{r,j}^2(x,\omega',\omega,E)d\omega'\leq M_1,\nonumber\\
&\int_{S'}\sigma_{r,j}^2(x,\omega,\omega',E) d\omega'\leq M_2,
\eea
for a.e. $(x,\omega,E)\in G\times S\times I$.
Finally, $K_{r,j}^3$ is of the form (as in section \ref{res-coll} above)
\[
(K_{r,j}^3\psi)(x,\omega,E)
=
\int_{I'}\int_{0}^{2\pi}
\sigma_{r,j}^3(x,E',E)
\psi(x,\gamma(E',E,\omega)(s),E')ds dE'
\]
where
$\gamma=\gamma(E',E,\omega):[0,2\pi]\to S$
be a parametrization of the curve
\[
\Gamma(E',E,\omega)=\{\omega'\in S\ |\ \omega'\cdot\omega-\mu(E',E)=0\}.
\]
${\sigma}_{r,j}^3:G\times I^2\to\R$ is a non-negative measurable function such that
\bea\label{ass-s-3}
&\int_{I'}{\sigma}_{r,j}^3(x,E',E)dE'\leq M_1, \nonumber\\
&\int_{I'}{\sigma}_{r,j}^3(x,E,E')dE'\leq M_2,
\eea
for a.e. $(x,E)\in G\times I$.

In a more concise way the components $K_{r,j}$ can be given by (see \cite{tervo18-up}, section 5)
\[
(K_{r,j}\psi)(x,\omega,E)=\sum_{k=1}^3\int_{I'}\int_{S'}\sigma_{kj}(x,\omega',\omega,E',E)\psi_k(x,\omega',E') d\rho^{kj}_{S}(\omega')d\rho^{kj}_I(E'),\quad j=1,2,3
\]
where we use in the relevant way the
pairs of measures $(\rho_I^{kj},\rho_{S}^{kj})=({\mc L}^1,\mu_{S})$ or
$(\rho_I^{kj},\rho_{S}^{kj})=({\mc L}^1,\mu_H^1)$ or
$(\rho_I^{kj},\rho_{S}^{kj})=(\mu_H^0,\mu_{S})$.
Here $\mc L^1$ ($=dE$) and $\mu_{S}$ ($=d\omega$) are the 1-dimensional Lebesgue measure on $I$ and the usual surface measure on $S$, respectively. $\mu_H^0$ and $\mu_H^1$ are the 0-dimensional and 1-dimensional Hadamard measures.
These notations serve only to compress the different
varieties of collision operators. Note that in this setting $\sigma_{kj}=\sigma_{r,j}^k$.
 
We recall from \cite{tervo18-up} (section 5) the following theorems 

\begin{theorem}\label{coup-bound}
Suppose that $\sigma_{kj}:G\times S^2\times I^2\to\R$ are non-negative measurable functions such that a.e. $(x,\omega,E)\in G\times S\times I$
\bea\label{bound-coupk}
&
\sum_{k=1}^3\int_{I'}\int_{S'}\sigma_{kj}(x,\omega',\omega,E',E)
d\rho_{S}^{kj}(\omega') d\rho_I^{kj}(E')\leq M_1
\nonumber\\
&
\sum_{k=1}^3\int_{I'}\int_{S'}\sigma_{jk}(x,\omega,\omega',E,E')
d\rho_{S}^{jk}(\omega') d\rho_I^{jk}(E')
\leq M_2
\eea
Then ${\bf K}_r: L^2(G\times S\times I)^3\to L^2(G\times S\times I)^3$ is a bounded operator and
\be
\n{{\bf K}_r}\leq 2\pi\sqrt{M_1M_2}. 
\ee
\end{theorem}

\begin{proof}   \end{proof}

\begin{theorem}\label{diss-for-coupled}
Suppose that for $j=1,2,3$ and a.e. $(x,\omega,E)\in G\times S\times I$
\bea\label{k3-n3-coupdiss}
&
\Sigma_j(x,\omega,E)-\sum_{k=1}^3\int_{ I'}\int_{S'}\sigma_{kj}(x,\omega',\omega,E',E)d\rho_{S}^{kj}(\omega') d\rho_I^{kj}(E')\geq c
\nonumber\\
&
\Sigma_j(x,\omega,E)-\sum_{k=1}^3\int_{I'}\int_{S'}\sigma_{jk}(x,\omega,\omega',E,E')d\rho_S^{jk}(\omega') d\rho_I^{jk}(E')\geq c 
\eea
Then for all $\psi\in L^2(G\times S\times I)^3$
\be\label{coerK-a}
\la ({\bf \Sigma}-{\bf K}_r)\psi,\psi\ra_{L^2(G\times S\times I)^3}\geq c\n{\psi}_{L^2(G\times S\times I)^3}^2.
\ee
\end{theorem}

\begin{proof}   \end{proof}

The Theorem \ref{coup-bound} is based to the Schur boundedness criterion (\cite{halmos}, p.22).
The conditions (\ref{ass-s-1}),
(\ref{ass-s-2}), (\ref{ass-s-3}) are equivalent to the assumptions (\ref{bound-coupk}). 

\begin{remark}
Notice that the coercivity (accretivity) estimate \eqref{coerK-a} is equivalent to the property that
for every $\lambda >0$,
and $\psi\in L^2(G\times S\times I)^3$,
\be\label{sda5a}
\n{ (\lambda I-(-{\bf \Sigma}+{\bf K}_r+cI))\psi}_{L^2(G\times S\times I)^3}\geq \lambda \n{\psi}_{L^2(G\times S\times I)^3},
\ee
which means that the operator $-{\bf \Sigma}+{\bf K}_r+cI:L^2(G\times S\times I)^3\to L^2(G\times S\times I)^3$ is {\it dissipative} (\cite[Section II.3.b]{engelnagel}, or \cite[Section 1.4]{pazy83}).
\end{remark}

For other data we impose the following criteria.
We assume that functions $a_j,\ b_j:\ol G\times I\to\R,\ j=2,3$
satisfy the following assumptions:

\be\label{trunc-18aa}
a_j\in L^\infty(G,W^{\infty,1}(I)), b_j\in L^\infty(G\times I),
\ee
\be\label{trunc-14aa}
-{\p {a_j}{E}}(x,E)
\geq q_{1,j}>0, \ {\rm a.e.},
\ee
\be\label{trunc-19aa}
-b_j(x,E)\geq q_{2,j}>0, \  {\rm a.e.},
\ee
\be\label{trunc-20aa}
-a_j(x,E_0)\geq q_{3,j}>0,\ 
-a_j(x,E_m)\geq q_{3,j}>0,  
\ {\rm a.e.} ,
\ee 
\be\label{trunc-20b}
|a_j(x,E)|\geq c_{0,j}>0,\ {\rm a.e.} .
\ee

As above, we apply the variational formulations to deduce existence of solutions. 
Recall that  the inner product in $L^2(G\times S\times I)^3$ is given by
\[
\la \phi,v\ra_{L^2(G\times S\times I)^3}
=\sum_{j=1}^3 \la \phi_j,v_j\ra_{L^2(G\times S\times I)},
\]
and analogously in other products of inner product spaces.
Integrating by parts and applying the Green's formula  we find (similarly as in section \ref{var-form}) that the bilinear form ${\bf B}(\cdot,\cdot):C^1(\ol G\times I,C^2(S))^3\times C^1(\ol G\times I,C^2(S))^3\to\R $ and the linear form ${\bf F}:
C^1(\ol G\times I,C^2(S))^3 \to\R$ corresponding to the problem 
(\ref{csda1a}), (\ref{csda1b}), (\ref{csda2}), (\ref{csda3})  are
(here we denote, as usual $\omega\cdot\nabla_x v:=(\omega\cdot\nabla_x v_2,\omega\cdot\nabla_x v_2,\omega\cdot\nabla_x v_3)$) and $\gamma_+(v):=(\gamma_+(v_1),\gamma_+(v_2),\gamma_+(v_3))$ and so on)
\bea\label{cosyst5}
{\bf B}(\psi,v)=&
-\sum_{j=2,3}\la \psi_j,a_j {\p {v_j}{E}}\ra_{L^2(G\times S\times I)} 
-\sum_{j=2,3}\la \psi_j,{\p {a_j}{E}}v_j \ra_{L^2(G\times S\times I)}\nonumber\\
&
-\la\psi,\omega\cdot\nabla_x v\ra_{L^2(G\times S\times I)^3}
-\sum_{j=2,3} \la \nabla_S\psi_j,b_j(x,E)
\nabla_Sv_j\ra_{L^2(G\times S\times I)}
\nonumber\\
&
+\la\psi,({\bf \Sigma}^*-{\bf K}_r^*) v\ra_{L^2(G\times S\times I)^3}
\nonumber\\
&
+\la \gamma_+(\psi), \gamma_+(v) \ra_{T^2(\Gamma_+)^3}
-\sum_{j=2,3} \la \psi_j(\cdot,\cdot,E_0),a_j(.,E_0) v_j(\cdot,\cdot,E_0)\ra_{L^2(G\times S)} ,
\eea
and
\be\label{cosyst6}
{\bf F}(v)
=
\la {\bf f},v\ra_{L^2(G\times S\times I)^3}
+\la {\bf g},\gamma_-(v)\ra_{T^2(\Gamma_-)^3}.
\ee

The appropriate solution spaces are (recall the definitions of $H_1,\ H_2$ and $\mc H,\ \widehat{\mc H}$ from sections \ref{b-f-spaces} and \ref{var-form}, respectively)
\be\label{eq:tilde_H}
&\mathbb{H}:=H_1\times \mc H\times \mc H, \nonumber \\
&\widehat{\mathbb{ H}}:=H_2\times \widehat{\mc H}\times \widehat{\mc H}.
\ee
The spaces $\mathbb{ H}$ and $\widehat {\mathbb{ H}}$ are Hilbert spaces equipped respectively with the inner products
\[
\la \psi,v\ra_{\mathbb{ H}}=\la\psi_1,v_1\ra_{H_1}+\sum_{j=2,3}\la\psi_j,v_j\ra_{\mc H}
\]
and
\[
\la\psi,v\ra_{\widehat{\mathbb{ H}}}
=\la \psi_1,v_1\ra_{H_2}+\sum_{j=2,3}\la\psi_j,v_j\ra_{\widehat{\mc H}}.
\]
The bilinear form ${\bf B}:C^1(\ol G\times I,C^2(S))^3\times C^1(\ol G\times I,C^2(S))^3\to\R$
has the following boundedness and coercivity properties.

\begin{theorem}\label{cosystth1}
Suppose that the assumptions (\ref{ass-S}),
(\ref{bound-coupk}), (\ref{k3-n3-coupdiss}), (\ref{trunc-18aa}), \ref{trunc-14aa}, (\ref{trunc-19aa}),   (\ref{trunc-20aa}) are valid.
Then there exists $M>0$ and $c'>0$ such that
\be\label{cosyst12}
|{\bf B}(\psi,v)|\leq M\n{\psi}_{{\mathbb{ H}}}\n{v}_{\widehat{\mathbb{ H}}}\quad \forall\psi, v\in C^1(\ol G\times I,C^2(S))^3,
\ee
and 
\be\label{cosyst13}
{\bf B}(\psi,\psi)\geq c'\n{\psi}_{\mathbb{ H}}^2\quad \forall\phi\in C^1(\ol G\times I,C^2(S))^3.
\ee

In addition, for all $v\in C^1(\ol G\times I,C^2(S))^3$,
\be\label{cosyst14}
|{\bf F}(v)|\leq 
\big(\n{\bf f}_{L^2(G\times S\times I)^3}+\n{\bf g}_{T^2(\Gamma_-)^3}\big)\n{v}_{\mathbb{ H}}.
\ee
\end{theorem}

\begin{proof}
Applying Theorems \ref{coup-bound} and \ref{diss-for-coupled}
the proof follows by minor modifications of the proof of Theorem 6.6 given in \cite{tervo19}. We omit details.
\end{proof}

Due to the above theorem, the bilinear form ${\bf B}(\cdot,\cdot)$ has a unique extension $\widetilde {\bf B}(\cdot,\cdot):\mathbb{ H}\times \widehat {\mathbb{H}}\to \R$ which satisfies
\be\label{cosyst15}
|\widetilde {\bf B}(\psi,v)|\leq M\n{\psi}_{\mathbb{ H}}\n{v}_{\widehat{ \mathbb{ H}}}\quad \forall\psi\in \mathbb{ H},\ v\in \widehat{\mathbb{H}},
\ee
and 
\be\label{cosyst16}
\widetilde {\bf B}(v,v)\geq c'\n{v}_{\mathbb{ H}}^2\quad \forall v\in \widehat{ \mathbb{ H}}.
\ee
Likewise, \eqref{cosyst14} implies that the linear form ${\bf F}$ has a unique bounded extension $\mathbb{ H}\to\R$, which we still denote by ${\bf F}$.
The variational equation corresponding to the problem 
(\ref{csda1a}), (\ref{csda1b}), (\ref{csda2}), (\ref{csda3}) is
\be\label{cosyst7}
\widetilde {\bf B}(\psi,v)={\bf F}(v)\quad \forall v\in \widehat {\mathbb{ H}}.
\ee

For the coupled transport system we have the following
 existence theorem which can be verified by similar ingredients as Theorem \ref{cor-csdath3-a} above.

\begin{theorem}\label{csdath3}
A. Suppose that the assumptions (\ref{ass-S}),
(\ref{bound-coupk}), (\ref{k3-n3-coupdiss}), (\ref{trunc-18aa}), \ref{trunc-14aa}, (\ref{trunc-19aa}), (\ref{trunc-20aa}) 
are valid. 
Let ${\bf f}\in L^2(G\times S\times I)^3$ and ${\bf g}\in T^2(\Gamma_-)^3$.
Then the  variational equation 
\be\label{csda40a}
\widetilde{\bf B}(\psi,v)={\bf F}(v)\quad \forall v\in \widehat{\mathbb {H}}
\ee
has a solution $\psi\in \mathbb{ H}$.
Furthermore, $\psi$ is a weak (distributional) solution of the system of equations
\begin{gather}
\omega\cdot\nabla_x\psi_1+\Sigma_1\psi_1-K_{r,1}\psi=f_1, \label{csda1a-1}\\
a_j(x,E){\p {\psi_j}{E}}
+b_j(x,E)\Delta_S\psi_j
+\omega\cdot\nabla_x\psi_j+\Sigma_{j}\psi_j-K_{r,j}\psi
=f_j,\quad j=2,3.\label{csda1b-2}
\end{gather}

B.  Suppose additionally to the assumptions of Part A that (\ref{trunc-20b})
holds. 
Then the system of equations (\ref{csda1a-1}), (\ref{csda1b-2}) together with the inflow boundary value and initial conditions
\be\label{ib-i-c}
\psi_{|\Gamma_-}={\bf g},\ \psi_j(.,.,E_m)=0 \ {\rm for}\ j=2,3
\ee
has a unique solution $\psi\in \mathbb{ H}$.

In addition,   the solution $\psi$ obeys the apriori  estimate
\be\label{csda40-syst}
\n{\psi}_{\mathbb{ H}}\leq C
\big(\n{{\bf f}}_{L^2(G\times S\times I)^3}+\n{{\bf g}}_{T^2(\Gamma_-)^3}\big).
\ee

\end{theorem}

\begin{proof}   \end{proof}

\begin{remark}\label{cg+}
Under the assumptions of Theorem \ref{csdath3} it also holds that
${\psi}_{|\Gamma_+}\in T^2(\Gamma_+)^3$ and  $\psi_j(.,.,E_0)\in L^2(G\times S),\ j=2,3$.
\end{remark}

\begin{remark}\label{coupth-re}
The existence of solutions for the problems (\ref{trunc-ex-3-b}) and (\ref{csda1a-1}), (\ref{csda1b-2}), (\ref{ib-i-c})
 can also be studied by applying $m$-dissipativity  and non-autonomous evolution operator theories (cf. \cite{tervo17}, \cite{tervo18}; especially Theorem 3.6 of \cite{tervo17}). The usual evolution operator theory entails that the data ${\bf f},\ {\bf g}$ is regular enough and that the compatibility condition like 
\be\label{comp}
g_j(.,.,E_m)=0 \ {\rm in}\ \Gamma_-\ {\rm for}\ j=2,3
\ee
holds. These methods may give more refined structure and regularity of the solution. Note that the compatibility condition (\ref{comp}) must  necessarily  be valid in order that the solution $\psi$ is continuous up to boundary.

\end{remark}

\subsection{Adjoint problem}\label{ad-prob}

Let
\[
&
P_1(x,\omega,E,D)\psi_1:=\omega\cdot\nabla_x\psi_1,\\
&
P_j(x,\omega,E,D)\psi_j:=
a_j(x,E){\p {\psi_j}{E}}
+b_j(x,E)\Delta_S\psi_j
+\omega\cdot\nabla_x\psi_j,\ j=2,3,
\]
\[
{\bf P}(x,\omega,E,D)\psi:=(P_1(x,\omega,E,D)\psi_1, P_2(x,\omega,E,D)\psi_2,P_3(x,\omega,E,D)\psi_3).
\]
Note that using this notation the system of equations (\ref{csda1a-1}), (\ref{csda1b-2})
is equivalent to
\[
{\bf P}(x,\omega,E,D)\psi+({\bf \Sigma}-{\bf K}_r)\psi={\bf f}.
\]
The formal adjoint ${\bf P}^*(x,\omega,E,D)$ of ${\bf P}(x,\omega,E,D)$ is
\[
{\bf P}^*(x,\omega,E,D)v:=(P_1^*(x,\omega,E,D)v_1, P_2^*(x,\omega,E,D)v_2,P_3^*(x,\omega,E,D)v_3)
\]
where
\[
&
P_1^*(x,\omega,E,D)v_1:=-\omega\cdot\nabla_xv_1,\\
&
P_j^*(x,\omega,E,D)v_j:=
-{\p {(a_jv_j)}{E}}
+b_j(x,E)\Delta_Sv_j
-\omega\cdot\nabla_xv_j,\ j=2,3
\]
for $v\in C_0^\infty(G\times S\times I^\circ)^3$.

The \emph{adjoint transport problem} to (\ref{csda1a-1}), (\ref{csda1b-2}),  (\ref{ib-i-c}) is formulated as 
\bea\label{ad-p-1}
&
-\omega\cdot\nabla_x\psi^*_1+\Sigma_1^*\psi^*_1-(K_{r}^*)_1\psi^*=f_1^*,\nonumber\\
&
-{\p {(a_j\psi^*_j)}{E}}
+b_j(x,E)\Delta_S\psi^*_j
-\omega\cdot\nabla_x\psi^*_j+\Sigma_{j}^*\psi^*_j-(K_{r}^*)_j\psi^*=f_j^*,\ j=2,3,
\eea
\be\label{ad-i-b}
\psi^*_{|\Gamma_+}={\bf g^*},\ \psi^*_j(.,.,E_0)=0 \ {\rm for}\ j=2,3
\ee
where ${\bf f}^*=(f_1^*,f_2^*,f_3^*)\in L^2(G\times S\times I)^3$ and ${\bf g}^*=(g_1^*,g_2^*,g_3^*)\in
T^2(\Gamma_+)^3$. Here $({ K}_{r}^*)_j$ are the components of the adjoint operator ${\bf K}_r^*:L^2(G\times S\times I)^3\to L^2(G\times S\times I)^3$ of ${\bf K}_r$. The components $({\bf K}_r^*)_j$ can be computed and 
\[
(({\bf K}_r^*)_j\psi^*)(x,\omega,E)=\sum_{k=1}^3\int_{I'}\int_{S'}\sigma_{jk}(x,\omega,\omega',E,E')\psi_k^*(x,\omega',E') d\rho^{jk}_{S}(\omega')d\rho^{jk}_I(E'),\quad j=1,2,3.
\]
The system of equations (\ref{ad-p-1})
is equivalent to
\[
{\bf P}^*(x,\omega,E,D)\psi^*+({\bf \Sigma}^*-{\bf K}^*_r)\psi^*={\bf f}^*.
\]

The 
bilinear   and linear forms corresponding to the adjoint problem 
(\ref{ad-p-1})
are respectively
\bea\label{ad-p-2}
&\widetilde{\bf B}^*(\psi^*,v)
=
\sum_{j=2,3}\la \psi_j^*,a_j {\p {v_j}{E}}\ra_{L^2(G\times S\times I)} 
+\la\psi^*,\omega\cdot\nabla_x v\ra_{L^2(G\times S\times I)^3}
\nonumber\\
&
-\sum_{j=2,3} \la \nabla_S\psi_j^*,b_j(x,E)
\nabla_Sv_j\ra_{L^2(G\times S\times I)}
+\la\psi^*,({\bf \Sigma}-{\bf K}_r) v\ra_{L^2(G\times S\times I)^3}
\nonumber\\
&
+\la \gamma_-(\psi^*), \gamma_-(v) \ra_{T^2(\Gamma_+)^3}
-\sum_{j=2,3} \la \psi_j^*(\cdot,\cdot,E_m),a_j(.,E_m) v_j(\cdot,\cdot,E_m)\ra_{L^2(G\times S)} ,
\eea
and
\be\label{ad-p-3}
{\bf F}^*(v)
=
\la {\bf f}^*,v\ra_{L^2(G\times S\times I)^3}
+\la {\bf g}^*,\gamma_+(v)\ra_{T^2(\Gamma_+)^3}.
\ee

For the adjoint problem, we have the following existence and uniqueness result
analogous to the above Theorem \ref{csdath3}.

\begin{theorem}\label{ad-th}
A. Suppose that the assumptions 
(\ref{ass-S}),
(\ref{bound-coupk}), (\ref{k3-n3-coupdiss}), (\ref{trunc-18aa}), \ref{trunc-14aa}, (\ref{trunc-19aa}), (\ref{trunc-20aa})
are valid. 
Let ${\bf f}^*\in L^2(G\times S\times I)^3$ and ${\bf g}^*\in T^2(\Gamma_+)^3$.
Then the  variational equation 
\be\label{ad-p-4}
\widetilde{\bf B}^*(\psi^*,v)={\bf F}^*(v)\quad \forall v\in \widehat{\mathbb {H}}
\ee
has a solution $\psi^*\in \mathbb{ H}$.
Furthermore, $\psi^*$ is a weak (distributional) solution of the system of equations (\ref{ad-p-1}).

B.  Suppose additionally to the assumptions of Part A that (\ref{trunc-20b})
holds  
Then the system of equations (\ref{ad-p-1}) together with the inflow boundary and initial conditions
\be\label{ad-i-c}
{\psi^*}_{|\Gamma_+}={\bf g}^*,\ \psi_j^*(.,.,E_0)=0 \ {\rm for}\ j=2,3
\ee
has a unique solution $\psi^*\in \mathbb{ H}$.
 
\end{theorem}

\begin{proof}  \end{proof}

From the operator theoretical point of view, the above existence result for adjoint problem is based on the fact that for a densely defined closed operator $A:X\to X$
in Hilbert space $X$ whose range $R(A^*)$ is closed in $X$, one has $R(A^*)=N(A)^{\perp}$ and $N(A^*)=R(A)^{\perp}$.
 
\begin{remark}\label{cg-}
Under the assumptions of Theorem \ref{ad-th} it also holds that
$\psi^*_{|\Gamma_-}\in T^2(\Gamma_-)^3$ and  $\psi_j^*(.,.,E_m)\in L^2(G\times S),\ j=2,3$.
\end{remark}

We need the next generalized Green's Theorem

\begin{theorem}\label{gen-green-a}
Suppose that the assumptions  
of Theorem \ref{csdath3} are valid.
Let $\psi$ and $\psi^*$ be the solutions of the problem (\ref{csda1a-1}), (\ref{csda1b-2}), (\ref{ib-i-c}) and its adjoint problem (\ref{ad-p-1}), (\ref{ad-i-b}). Then
\bea\label{g-g-1}
&
\la {\bf P}(x,\omega,E,D)\psi,\psi^*\ra_{L^2(G\times S\times I)^3}
-
\la \psi,{\bf P}^*(x,\omega,E,D)\psi^*\ra_{L^2(G\times S\times I)^3} 
\nonumber\\
&
=
\sum_{j=1}^3\int_{\Gamma}\psi_j\psi_j^*(\omega\cdot\nu) dx d\sigma dE
\nonumber\\
&
+
\sum_{j=2,3} \la \psi_j(\cdot,\cdot,E_m),a_j(.,E_m) \psi^*_j(\cdot,\cdot,E_m)\ra_{L^2(G\times S)}\nonumber\\
&
-\sum_{j=2,3} \la \psi_j(\cdot,\cdot,E_0),a_j(.,E_0) \psi^*_j(\cdot,\cdot,E_0)\ra_{L^2(G\times S)} .
\eea
\end{theorem}

\begin{proof}
By Remark \ref{cg+} $\gamma_+(\psi)\in T^2(\Gamma_+)^3$ and  $\psi_j(\cdot,\cdot,E_0)\in L^2(G\times S),\ j=2,3$ and by By Remark \ref{cg-} $\gamma_-(\psi^*)\in T^2(\Gamma_-)^3$ and  $\psi^*_j(\cdot,\cdot,E_m)\in L^2(G\times S),\ j=2,3$. Note also that $ {\bf P}(x,\omega,E,D)\psi,\ {\bf P}^*(x,\omega,E,D)\psi^*\in L^2(G\times S\times I)^3$. The proof follows by limiting techniques from the standard Green formula for smooth functions. We omit details.
\end{proof}

The solutions $\psi$ and $\psi^*$ are related in the following way 

\begin{lemma}\label{le-sol-id}
Suppose that the assumptions (\ref{ass-S}),
(\ref{bound-coupk}), (\ref{k3-n3-coupdiss}), (\ref{trunc-18aa}), \ref{trunc-14aa}, (\ref{trunc-19aa}), (\ref{trunc-20aa}), (\ref{trunc-20b})
(of Theorem \ref{csdath3}) are valid.
Let $\psi$ and $\psi^*$ be the solutions of the problem (\ref{csda1a-1}), (\ref{csda1b-2}), (\ref{ib-i-c}) and its adjoint problem (\ref{ad-p-1}), (\ref{ad-i-b}). 
Then  the expressions $\widetilde{\bf B}(\psi,\psi^*)$ and $\widetilde{\bf B}^*(\psi^*,\psi)$ are well-defined and
\be\label{sol-id-a}
\widetilde{\bf B}(\psi,\psi^*)={\bf F}(\psi^*),
\ee
\be\label{sol-id-b}
\widetilde{\bf B}^*(\psi^*,\psi)={\bf F}^*(\psi)
\ee
and
\be\label{sol-id}
\widetilde{\bf B}(\psi,\psi^*)=\widetilde{\bf B}^*(\psi^*,\psi).
\ee
\end{lemma}

\begin{proof}
 
A.
By partial integration and by the Green's Theorem for $v,\ v^*\in C^1(\ol G\times I\times S)$.
\[
&
-\la {\p {(a_jv_j^*)}{E}},v_j\ra_{L^2(G\times S\times I)}
=\la a_j{\p {v_j}{E}},v_j^*\ra_{L^2(G\times S\times I)}\\
&
-\big(\la v_j(\cdot,\cdot,E_m),a_j(.,E_m) v_j^*(\cdot,\cdot,E_m)\ra_{L^2(G\times S)}-
\la v_j(\cdot,\cdot,E_0),a_j(.,E_0) v_j^*(\cdot,\cdot,E_0)\ra_{L^2(G\times S)}
\big)
\]
and
\[
&
-\la v,\omega\cdot\nabla_x v^*\ra_{L^2(G\times S\times I)^3}
=
\la \omega\cdot\nabla_x v,v^*\ra_{L^2(G\times S\times I)^3}
\\
&
-
\big(\la \gamma_+(v),\gamma_+(v^*)\ra_{T^2(\Gamma_+)^3}
-\la \gamma_-(v),\gamma_-(v^*)\ra_{T^2(\Gamma_+)^3}\big).
\]
Hence for $v,\ v^*\in C^1(\ol G\times I\times S)$
\bea\label{sol-id-1}
&{\bf B}(v,v^*)=
-\sum_{j=2,3}\la v_j,{\p {(a_jv^*_j)}{E}}\ra_{L^2(G\times S\times I)} 
\nonumber\\
&
-\la v,\omega\cdot\nabla_x v^*\ra_{L^2(G\times S\times I)^3}
-\sum_{j=2,3} \la \nabla_S v_j,b_j(x,E)
\nabla_S v^*_j\ra_{L^2(G\times S\times I)}
\nonumber\\
&
+\la v,({\bf \Sigma}^*-{\bf K}_r^*) v^*\ra_{L^2(G\times S\times I)^3}
\nonumber\\
&
+\la \gamma_+(v), \gamma_+(v^*) \ra_{T^2(\Gamma_+)^3}
-\sum_{j=2,3} \la v_j(\cdot,\cdot,E_0),a_j(.,E_0) v^*_j(\cdot,\cdot,E_0)\ra_{L^2(G\times S)} \nonumber\\
&
=
\sum_{j=2,3}
\la a_j{\p {v_j}{E}},v_j^*\ra_{L^2(G\times S\times I)}
\nonumber\\
&
-\sum_{j=2,3}\big(\la v_j(\cdot,\cdot,E_m),a_j(.,E_m) v_j^*(\cdot,\cdot,E_m)\ra_{L^2(G\times S)}-
\la v_j(\cdot,\cdot,E_0),a_j(.,E_0) v_j^*(\cdot,\cdot,E_0)\ra_{L^2(G\times S)}
\big)\nonumber\\
&
+
\la \omega\cdot\nabla_x v,v^*\ra_{L^2(G\times S\times I)^3}
\nonumber\\
&
-
\big(\la \gamma_+(v),\gamma_+(v^*)\ra_{T^2(\Gamma_+)^3}
-\la \gamma_-(v),\gamma_-(v^*)\ra_{T^2(\Gamma_+)^3}\big)\nonumber\\
&
-\sum_{j=2,3} \la \nabla_S v_j,b_j(x,E)
\nabla_S v^*_j\ra_{L^2(G\times S\times I)}
\nonumber\\
&
+\la v,({\bf \Sigma}^*-{\bf K}_r^*) v^*\ra_{L^2(G\times S\times I)^3}
\nonumber\\
&
+\la \gamma_+(v), \gamma_+(v^*) \ra_{T^2(\Gamma_+)^3}
-\sum_{j=2,3} \la v_j(\cdot,\cdot,E_0),a_j(.,E_0) v^*_j(\cdot,\cdot,E_0)\ra_{L^2(G\times S)} \nonumber\\
&
=
\sum_{j=2,3}
\la a_j{\p {v_j}{E}},v_j^*\ra_{L^2(G\times S\times I)}
+
\la \omega\cdot\nabla_x v,v^*\ra_{L^2(G\times S\times I)^3}
\nonumber\\
&
-\sum_{j=2,3} \la \nabla_S v_j,b_j(x,E)
\nabla_S v^*_j\ra_{L^2(G\times S\times I)}
\nonumber\\
&
+\la ({\bf \Sigma}-{\bf K}_r) v,v^*\ra_{L^2(G\times S\times I)^3}
\nonumber\\
&
+\la \gamma_-(v), \gamma_-(v^*) \ra_{T^2(\Gamma_+)^3}
-\sum_{j=2,3} \la a_j(.,E_m) v_j(\cdot,\cdot,E_m), v^*_j(\cdot,\cdot,E_m)\ra_{L^2(G\times S)}
={\bf B}^*(v^*,v).
\eea

For all $v,\ v^*\in C^1(\ol G\times I,C^2( S))$
\[
&
-\la \nabla_S v_j,b_j(x,E)
\nabla_S v^*_j\ra_{L^2(G\times S\times I)}
\\
&
=
\la  v_j,b_j(x,E)
\Delta_S v^*_j\ra_{L^2(G\times S\times I)}
=
\la b_j(x,E)
\Delta_S v_j,v_j^*\ra_{L^2(G\times S\times I)}
\]
and so by
using the above notations we obtain by (\ref{sol-id-1}) 
\bea\label{sol-id-1a}
&
{\bf B}(v,v^*)=
\la v,{\bf P}^*(x,\omega,E,D)v^*\ra_{L^2(G\times S\times I)^3} 
+\la v,({\bf \Sigma}^*-{\bf K}_r^*) v^*\ra_{L^2(G\times S\times I)^3}
\nonumber\\
&
+\la \gamma_+(v), \gamma_+(v^*) \ra_{T^2(\Gamma_+)^3}
-\sum_{j=2,3} \la v_j(\cdot,\cdot,E_0),a_j(.,E_0) v^*_j(\cdot,\cdot,E_0)\ra_{L^2(G\times S)} \nonumber\\
&
=
\la {\bf P}(x,\omega,E,D)v,v^*\ra_{L^2(G\times S\times I)^3}
+\la ({\bf \Sigma}-{\bf K}_r) v,v^*\ra_{L^2(G\times S\times I)^3}
\nonumber\\
&
+\la \gamma_-(v), \gamma_-(v^*) \ra_{T^2(\Gamma_-)^3}
-\sum_{j=2,3} \la a_j(.,E_m)v_j(\cdot,\cdot,E_m), v^*_j(\cdot,\cdot,E_m)\ra_{L^2(G\times S)}\nonumber\\
&
={\bf B}^*(v^*,v).
\eea

B.
By Theorem \ref{gen-green-a} we find that (\ref{sol-id-1a}) holds also for 
$v=\psi,\ v^*=\psi^*$ and so (\ref{sol-id}) is valid. Furthermore, recalling  that ${\bf P}(x,\omega,E,D)\psi+
(\Sigma-K_r) \psi={\bf f}$, $\gamma_-(\psi)={\bf g}$ we have
\bea\label{soli-2}
&
{\bf B}(\psi,\psi^*)=
\la {\bf P}(x,\omega,E,D)\psi,\psi^*\ra_{L^2(G\times S\times I)^3}
+\la ({\bf \Sigma}-{\bf K}_r) \psi,\psi^*\ra_{L^2(G\times S\times I)^3}
\nonumber\\
&
+\la \gamma_-(\psi), \gamma_-(\psi^*) \ra_{T^2(\Gamma_+)^3}
-\sum_{j=2,3} \la \psi_j(\cdot,\cdot,E_m),a_j(.,E_m) \psi^*_j(\cdot,\cdot,E_m)\ra_{L^2(G\times S)}\nonumber\\
&
=
\la {\bf f},\psi^*\ra_{L^2(G\times S\times I)^3}
+\la {\bf g}, \gamma_-(\psi^*) \ra_{T^2(\Gamma_-)^3}
={\bf F}(\psi^*)
\eea
and so (\ref{sol-id-a}) holds. The claim (\ref{sol-id-b}) is similarly proved which completes the proof.

\end{proof}


\section{Inverse radiation treatment problem}\label{irtpre}

\subsection{Radiation Treatment Planning}\label{RTP}

Radiation therapy aims to generate dose distributions in such a way that the desired dose conforms to the target volume, while the healthy tissue and especially the so-called critical organs achieve as low dose as possible. 
Dose can be delivered externally (external therapy) or internally (internal therapy or brachytherapy).
The determination of (optimal) incoming external particle fluxes through the patches of patient surface or internal sources located inside the patient tissue is the basic task in treatment planning known as the {\it inverse radiation treatment planning.} We recall some details from \cite{tervo17-up}.

The patient domain $G\subset\R^3$ consists of the tumor volume ${\s T}$,
the critical organ region ${\s C}$ and the normal tissue region ${\s N}$,
as a mutually disjoint union $G={\s T}\cup {\s C}\cup {\s N}$. We assume that the sets
${\s T},\ {\s C},\ {\s N}$ are Lebesgue measurable.
The tumor volume, that is the target, includes the tumor and some safety margin around it.
Critical organs and normal tissue are build up of healthy tissue,
and must be conserved during the treatment as well as possible.
In the sequel we only deal with the inverse treatment planning problem in the context
of the stationary Boltzmann transport equation (i.e. we omit time dependency here which could be treated analogously after relevant modifications).

In radiation therapy the dose absorbed from the particle field $\psi=(\psi_1,\psi_2,\psi_3)$ is defined by the functional
\be\label{inv-1}
D(x)=(D\psi)(x):=\sum_{j=1}^3\int_{S\times I}\varsigma_j (x,E)\psi_j(x,\omega,E) d\omega dE,
\ee
where $\varsigma_j\in L^\infty(G\times I)$, $\varsigma_j\geq 0$ are the so-called \emph{(total) stopping powers}.
The component fields of $\psi$, relevant to photon and electron radiation therapy, are $\psi_1={}$photons, $\psi_2={}$electrons and $\psi_3={}$positrons.
We find that $D:L^2(G\times S\times I)^3\to L^2(G)$ is a bounded linear operator and  its adjoint operator $D^*:L^2(G)\to L^2(G\times S\times I)^3$ is simply a multiplication type operator,
\[
D^*d=(\varsigma_1,\varsigma_2,\varsigma_3)d,\quad {\rm for}\ d\in L^2(G).
\]
We describe below shortly an optimization problem related 
to inverse radiation treatment planning. 

In optimization
the dose $D=D(x)$ is computed from  (\ref{inv-1})
where $\psi=(\psi_1,\psi_2,\psi_3)\in L^2(G\times S\times I)^3$ satisfies
the \emph{coupled system of  transport equations},
\begin{gather}
\omega\cdot\nabla_x\psi_1+\Sigma_1\psi_1-K_{r,1}\psi=f_1, \label{inv-2}\\
a_j(x,E){\p {\psi_j}{E}}
+b_j(x,E)\Delta_S\psi_j
+\omega\cdot\nabla_x\psi_j+\Sigma_{j}\psi_j-K_{r,j}\psi
=f_j,\quad j=2,3.\label{inv-3}
\end{gather}
\be\label{inv-4}
\psi_{|\Gamma_-}={\bf g},\ \psi_j(.,.,E_m)=0,\ j=2,3.
\ee

We formulate the equations (\ref{inv-2}), (\ref{inv-3}), (\ref{inv-4}) in the more concise way. 
Let
\[
&
{\bf A}\psi:=\Big(\omega\cdot\nabla_x\psi_1+\Sigma_1\psi_1-K_{r,1}\psi,
\\
&
a_2(x,E){\p {\psi_2}{E}}
+b_2(x,E)\Delta_S\psi_2
+\omega\cdot\nabla_x\psi_2+\Sigma_{2}\psi_2-K_{r,2}\psi,
\\
&
a_3(x,E){\p {\psi_3}{E}}
+b_3(x,E)\Delta_S\psi_3
+\omega\cdot\nabla_x\psi_3+\Sigma_{3}\psi_2-K_{r,3}\psi\Big).
\]
Moreover, 
let ${\bf T}:=({\bf A},\gamma_-):L^2(G\times S\times I)^3\to L^2(G\times S\times I)^3\times T(\Gamma_-)^3$ be the unbounded linear operator with the domain
\[
D({\bf T})&:=\{\psi\in \mathbb{H}|\ {\bf T}\psi=({\bf A}\psi,\gamma_-(\psi))=({\bf f},{\bf g})\in 
L^2(G\times S\times I)^3\times T(\Gamma_-)^3
,\ \\
&
\psi_2(.,.,E_m)=\psi_3(.,.,E_m)=0\}.
\]
Then $\psi$ is a solution of the problem (\ref{inv-2}), (\ref{inv-3}), (\ref{inv-4}) if and only if
\be\label{inv-6}
{\bf T}\psi=({\bf f},{\bf g})
\ee
and under the asumptions of Theorem \ref{csdath3} ${\bf T}$ is invertible and so
\be\label{inv-7}
\psi={\bf T}^{-1}({\bf f},{\bf g})=:\psi({\bf f},{\bf g}).
\ee
The generated dose is then
\[
D=D(x)=\big(D({\bf T}^{-1}({\bf f},{\bf g}))\big)(x),\quad x\in G.
\]
We denote 
\[
{\s D}({\bf f},{\bf g}):=D(\psi({\bf f},{\bf g})).
\] 
Then we have

\begin{lemma}\label{dc}
The dose operator ${\s D}:L^2(G\times S\times I)^3\times T^2(\Gamma_-)^3\to L^2(G)$ is linear and bounded, i.e.
\[
\n{{\s D}({\bf f},{\bf g})}_{L^2(G)}\leq C\ (\n{{\bf f}}_{L^2(G\times S\times I)^3}+\n{{\bf g}}_{ T^2(\Gamma_-)^3}).
\]
\end{lemma}

\begin{proof}

We see that ${\s D}({\bf f},{\bf g})=D({\bf T}^{-1}({\bf f},{\bf g}))$ and so 
${\s D}$ is linear.
Since $D:L^2(G\times S\times I)^3\to L^2(G)$ is bounded we have by
(\ref{csda40-syst})
\[
&
\n{{\s D}({\bf f},{\bf g})}_{L^2(G)}=\n{{D}(\psi({\bf f},{\bf g}))}_{L^2(G)}
\leq 
\n{D}\n{(\psi({\bf f},{\bf g}))}_{L^2(G\times S\times I)^3}
\\
&
\leq 
\n{D}\n{(\psi({\bf f},{\bf g}))}_{\mathbb{H}}
\leq 
\n{D}C
\big(\n{{\bf f}}_{L^2(G\times S\times I)^3}+\n{{\bf g}}_{T^2(\Gamma_-)^3}\big)
\]
which completes the proof.

\end{proof}

Commonly used {\it physical criteria} for the dose
are the following ones. We demand that
\begin{align}
&D(x)=D_0,\ x\in {\s T},\label{tc-a}\\
&D(x)\leq D_C,\ x\in {\s C},\label{cc}\\  
&D(x)\leq D_N,\ x\in {\s N},\label{nc}
\end{align}
where $D_0$ is the prescribed (usually uniform) dose in the target ${\s T}$ and where $D_C$ and $D_N$ are the allowed upper bounds (usually uniforms as well) in the critical organ ${\s C}$ and normal tissue ${\s N}$ regions, respectively. Instead of (\ref{tc-a}) one may require for more flexibly (when considering the so-called feasible solutions) that only
\bea
d_T\leq D(x)\leq D_T,\ x\in {\s T},\label{tc1}
\eea
where $D_T$ and $d_T$ are upper and lower bounds for dose in the target.

In addition to the above requirements in modern planning, one imposes so-called {\it dose volume constraints} for the dose distribution, especially for the critical organ region (but also for some other tissue region's similar dose volume constraints may be considered). Dose volume constraint demands that the dose $D(x)$ cannot be greater than some prescribed dose level, say $d_C$, in a prescribed volume fraction $v_C$
of ${\s C}$ . 
This can be expressed as follows
\bea\label{dv}
{{m(\{x\in {\s C}|\ D(x)\geq d_C\})}\over{m({\s C})}}\leq v_C
\eea
where $m$ is the 3-dimensional Lebesgue measure.
Clearly the dose volume constraint is equivalent to
\be\label{dva}
{1\over{m({\s C})}}\int_{\s C}H(D(x)-d_C) dx\leq v_C
\ee
where $H$ is the Heaviside function. Note that the integral in (\ref{dva}) exists.
In practice $H$ can be replaced  with a smooth 
function $H_\epsilon$ which approximates it to some reasonable level of accuracy.

\subsection{Object Function and the Optimization Problem}\label{objectf}
   
Our aim is that the above requirements \eqref{tc-a}, \eqref{cc}, \eqref{nc}, \eqref{dv} for the dose distribution are valid as well as possible.
For that purpose, we define the object (cost) function
\bea\label{of}
J({\bf f},{\bf g})=c_{\s T} J_{\s T}({\bf f},{\bf g})+c_{\s C} J_{\s C}({\bf f},{\bf g})+c_{\s N} J_{\s N}({\bf f},{\bf g})+c_{\rm DV} J_{\rm DV}({\bf f},{\bf g}),
\eea 
where
\[
J_{\s T}({\bf f},{\bf g})&=\n{D_0-{\s D}({\bf f},{\bf g})}_{L^2({\s T})}^2,\\
J_{\s C}({\bf f},{\bf g})&=\n{(D_C-{\s D}({\bf f},{\bf g}))_-}_{L^2({\s C})}^2,\\
J_{\s N}({\bf f},{\bf g})&=\n{(D_N-{\s D}({\bf f},{\bf g}))_-}_{L^2({\s N})}^2,\\
J_{\rm DV}({\bf f},{\bf g})&=
\Big(\big(v_C-{1\over{m({\s C})}}\int_{\s C}H({\s D}({\bf f},{\bf g})(x)-d_C) dx\big)_-\Big)^2,
\]
and where $c_{\s T}, c_{\s C}, c_{\s N}, c_{\rm DV}$ are non-negative weights
with which one controls the different 
priorities in the optimization. Here $a_-$ denotes the \emph{negative part} 
of $a\in\R$,
\[
a_-=-\min\{a,0\}=\frac{1}{2}(|a|-a).
\]

In external radiotherapy ${\bf f}=0$ and in internal radiotherapy ${\bf g}=0$
(from the mathematical point of view both can, of course, be non-zero and the below considerations apply after relevant modifications for this more general case). Thus the corresponding object functions in usual practice are
\[
J_{\rm ex}({\bf g}):=&J(0,{\bf g})\ \ {\rm (external\ radiotherapy})\ {\rm and}\\
J_{\rm in}({\bf f}):=&J({\bf f},0)\ \ {\rm (internal\ radiotherapy)}.
\]

The {\it admissible sets}
for the optimal control problems  are respectively
\[
U_{\rm ad}=\{{\bf g}\in T^2(\Gamma_-)^3\ |\ {\bf g}\geq 0\}\ ({\rm external\ radiotherapy})
\]
and
\[
U_{\rm ad}'=\{{\bf f}\in L^2(G\times S\times I))^3\ |\ {\bf f}\geq 0\}\ ({\rm internal\ radiotherapy}).
\]  
They both are convex sets (cones) of the ambient spaces. If (in practical optimization) one wants to take the whole ambient space as an admissible set that is, $U_{\rm ad}=T^2(\Gamma_-)^3$ or
respectively $U_{\rm ad}'=L^2(G\times S\times I)^3$
one must add to the object function the penalty term
\[
+c_{\rm ad} J_{\rm ad}({\bf g}),
\ {\rm where},\ J_{\rm ad}({\bf g})=\n{{\bf g}_-}^2_{T^2(\Gamma_-)^3}
\ {(\rm external\ radiotherapy)}
\]
and
\[
+c_{\rm ad'} J_{\rm ad'}({\bf f}),
\ {\rm where}\ J_{\rm ad'}({\bf f})=\n{{\bf f}_-}^2_{L^2(G\times S\times I)^3}
\ {\rm (internal\ radiotherapy)}.
\]
These take care of the non-negativity of the incoming flux or source, respectively.
In theory as well as in practice, it is also reasonable to add \emph{stabilizing cost terms}
correspondingly
\be\label{ofsc1}
+c_{\rm sc} J_{\rm sc}({\bf g})\ {\rm where}\ J_{\rm sc}({\bf g})=\n{\psi(0,{\bf g})}^2_{L^2(G\times S\times I)^3}\ {\rm or}\ J_{\rm sc}({\bf g})=\n{{\bf g}}^2_{T^2(\Gamma_-)^3},
\ee
or
\be\label{ofsc2}
+c_{\rm sc'} J_{\rm sc'}({\bf f})\ {\rm where}\ J_{\rm sc'}({\bf f})=\n{\psi({\bf f},0)}^2_{L^2(G\times S\times I)^3}
\ {\rm or}\ J_{\rm sc'}({\bf f})=\n{{\bf f}}^2_{L^2(G\times S\times I)^3}.
\ee

As a conclusion, we have for the external therapy the object function
\bea\label{fof}
J_{\rm ex}({\bf g})=c_{\s T} J_{\s T}(0,{\bf g})+c_{\s C} J_{\s C}(0,{\bf g})+c_{\s N}J_{\s N}(0,{\bf g})
+c_{\rm DV} J_{\rm DV}(0,{\bf g})+c_{\rm sc}J_{\rm sc}({\bf g}),
\eea
when $U_{\rm ad}=\{{\bf g}\in T^2(\Gamma_-)^3|\ {\bf g}\geq 0\}$ or
\bea\label{fofa}
J_{\rm ex}({\bf g})=c_{\s T}J_{\s T}(0,{\bf g})+c_{\s C}J_{\s C}(0,{\bf g})+c_{\s N}J_{\s N}(0,{\bf g})
+c_{\rm DV}J_{\rm DV}(0,{\bf g})+c_{\rm ad} J_{\rm ad}({\bf g})+c_{\rm sc}J_{\rm sc}({\bf g}).
\eea
when $U_{\rm ad}=T^2(\Gamma_-)^3$.
The object function $J_{\rm in}({\bf f})$ for the internal therapy is formulated analogously.

With these concepts the overall \emph{optimal control initial inflow boundary value problem} can be stated as
: Find the {\it global minimum}
\bea\label{taske}
\min\{J_{\rm ex}({\bf g})\ |\ {\bf g}\in U_{\rm ad}\}\quad {(\rm external\ radiotherapy)},
\eea
or
\bea\label{taski}
\min\{J_{\rm in}({\bf f})\ |\ {\bf f}\in U'_{\rm ad}\}\quad {\rm (internal\ radiotherapy)},
\eea
where $U_{\rm ad}$ or $U'_{\rm ad}$, respectively, is chosen in the way explained above.

In \cite{tervo17-up}, section 8 we showed some features (such as convexity and Lipschitz continuity properties) of the above object functions. The search of global minimum requires usually \emph{global optimization methods} (\cite{pinter}, \cite{pinter14}). Global optimization algorithms in turn require an \emph{initial point}. In the following section we shall give one option to calculate a potential initial point.

\subsection{Computation of Initial Solutions}\label{cis}

One possibility to help the optimization process is to compute the initial solution for actual (global) optimization scheme as accurately and rapidly as possible.
We suggest the following approach, and consider its details mainly
in the case of the external radiotherapy.
The computations for internal radiotherapy (formulated below) are analogous
and thus omitted.
The formulas related to the
optimal control system 
are written below by using the relevant \emph{variational equations}.
The approach has some obvious benefits.
For example, the finite element method (FEM) schemes can be naturally implemented applying this formulation.

Let ${\bf f}\in L^2(G\times S\times I)^3$ and ${\bf g}\in T^2(\Gamma_-)^3$. Recall  that
under the assumptions of Theorem \ref{csdath3} $\psi$ is a solution of the problem
\begin{gather}
\omega\cdot\nabla_x\psi_1+\Sigma_1\psi_1-K_{r,1}\psi=f_1, \label{inv-2a}\\
a_j(x,E){\p {\psi_j}{E}}
+b_j(x,E)\Delta_S\psi_j
+\omega\cdot\nabla_x\psi_j+\Sigma_{j}\psi_j-K_{r,j}\psi
=f_j,\quad j=2,3.\label{inv-3a}
\end{gather}
\be\label{inv-4a}
\psi_{|\Gamma_-}={\bf g},\ \psi_j(.,.,E_m)=0,\ j=2,3.
\ee 
if and only if
\be\label{va1}
\widetilde{\bf B}(\psi,v)={\bf F}(v),\quad \forall v\in \widehat{\mathbb{H}}
\ee
where $\widetilde{\bf B}(\cdot,\cdot)$ is the bilinear form given by 
\bea\label{va2}
\widetilde{\bf B}(\psi,v)=&
-\sum_{j=2,3}\la \psi_j,a_j {\p {v_j}{E}}\ra_{L^2(G\times S\times I)} 
-\sum_{j=2,3}\la \psi_j,{\p {a_j}{E}}v_j \ra_{L^2(G\times S\times I)}\nonumber\\
&
\nonumber\\
&
-\la\psi,\omega\cdot\nabla_x v\ra_{L^2(G\times S\times I)^3}
-\sum_{j=2,3} \la \nabla_S\psi_j,b_j(x,E)
\nabla_Sv_j\ra_{L^2(G\times S\times I)}
\nonumber\\
&
+\la\psi,({\bf \Sigma}^*-{\bf K}_r^*) v\ra_{L^2(G\times S\times I)^3}
\nonumber\\
&
+\la \gamma_+(\psi), \gamma_+(v) \ra_{T^2(\Gamma_+)^3}
-\sum_{j=2,3} \la \psi_j(\cdot,\cdot,E_0),a_j(.,E_0) v_j(\cdot,\cdot,E_0)\ra_{L^2(G\times S)} .
\eea
Here ${\bf F}$
is the linear form
\be\label{va3}
{\bf F}(v)=\la {\bf f},v\ra_{L^2(G\times S\times I)^3}
+\la {\bf g},v\ra_{T^2(\Gamma_-)^3},\ v\in {\mathbb{H}}.
\ee
and  $\Sigma^* = \Sigma$ and $K_r^*$ is the adjoint of $K_r$.
Note that the bilinear form $\widetilde{\bf B}$ is not symmetric.
The existence result of the solution to the problem (\ref{va1}) or equivalently to the problem (\ref{inv-2a}), (\ref{inv-3a}), (\ref{inv-4a}) has been given above in Theorem \ref{csdath3}.

In external radiotherapy ${\bf f}=0$ and so we have $\psi=\psi(0,{\bf g})$ and ${\s D}(0,{\bf g})=D(\psi(0,{\bf g}))$.
We shall denote more shortly $\psi({\bf g})=\psi(0,{\bf g})$ and ${\s D}({\bf g})={\s D}(0,{\bf g})$.
Suppose that $d_{\s T}\in L^2({\s T}),\ d_{\s C}\in L^2({\s C}),\ d_{\s N}\in L^2({\s N})$ are
some given dose distributions (for example, they may be constants). We define an object function
\be\label{va9}
J({\bf g})=&c_{\s T}\n{ d_{\s T}-{\s D}({\bf g})}_{L^2({\s T})}^2+c_{\s C}\n{ d_{\s C}-{\s D}({\bf g})}_{L^2({\s C})}^2 \\
&+c_{\s N}\n{d_{\s N}-{\s D}({\bf g})}_{L^2({\s N})}^2+c_{\rm sc}\n{{\bf g}}_{T^2(\Gamma_-)^3}^2, \nonumber
\ee
with strictly positive constants $c_{\s T}, c_{\s C}, c_{\s N}, c_{\rm sc}>0$.
In practice, only one type of particles are inflowing simultaneously (usually photons or electrons) and so only one component of ${\bf g}$ is non-zero but we can formulate a more general result which includes this
more realistic situation. We denote 
\[
{\bf F}(v)=({\bf F}{\bf g})(v)
:=\sum_{j=1}^3\int_{\partial G\times S\times I}(\omega\cdot \nu)_-g_jv_j d\sigma d\omega dE
=\la {\bf g},\gamma_-(v)\ra_{T^2(\Gamma_-)^3}
\]
and, as before,
\[
U_{\rm ad}=\{{\bf g}\in T^2(\Gamma_-)^3\ |\ {\bf g}\geq 0\}.
\] 

We shall write 
\[
a_+=\max\{0,a\}={1\over 2}(a+|a|)
\]
for the positive part of $a\in\R$, and $a_+=((a_1)_+,(a_2)_+,(a_3)_+)$
when $a=(a_1,a_2,a_3)\in\R^3$.
As in \cite{tervo17-up} (Theorem 8.7) we have the following optimality result.

\begin{theorem}\label{thinit}
Suppose that the assumptions 
(\ref{ass-S}),
(\ref{bound-coupk}), (\ref{k3-n3-coupdiss}), (\ref{trunc-18aa}), \ref{trunc-14aa}, (\ref{trunc-19aa}), (\ref{trunc-20aa}), (\ref{trunc-20b})
(of Theorem \ref{csdath3})  
are satisfied. Then the minimum $\min_{{\bf g}\in U_{\rm ad}}{J}({\bf g})$ exists and it is given by
\be\label{v9a}
\ol {\bf g}={1\over {c_{\rm sc}}}(\gamma_-(\psi^*))_+,
\ee
where the pair $(\psi,\psi^*)\in \mathbb{H}\times \mathbb{H}$ is the unique solution of the \emph{coupled (non-linear) system} of variational equations (recall the dose operator $D$ from (\ref{inv-1}))
\bea\label{v10}
\widetilde{\bf B}^*(\psi^*,v)&+ c_{\s T}\la D\psi,Dv\ra_{ L^2({\s T})}+c_{\s C}\la D\psi,Dv\ra_{ L^2({\s C})}+c_{\s N}\la D\psi,Dv\ra_{ L^2({\s N})}\nonumber\\
&=c_{\s T}\la d_{\s T},Dv\ra_{ L^2({\s T})}+c_{\s C}\la d_{\s C},Dv\ra_{ L^2({\s C})}+c_{\s N}\la d_{\s N},Dv\ra_{ L^2({\s N})} \\
\widetilde{\bf B}(\psi,v)&=\la {1\over {c_{\rm sc}}}(\gamma_-(\psi^*))_+,\gamma_-(v)\ra_{T^2(\Gamma_-)^3},
\nonumber
\eea
for all  $v\in \widehat{\mathbb{H}}$.
\end{theorem}

\begin{proof}
By Lemma \ref{dc}  the dose operator ${\s D}:T^2(\Gamma_-)^3\to L^2(G)$ is a bounded linear operator. Hence the object function $J:T^2(\Gamma_-)^3\to \R$ is differentiable and strictly convex (see e.g. the proof of Theorem 8.3  given in \cite{tervo17-up}; note that the last term of $J$ makes it strictly convex).
Furthermore, $J$ is bounded from below (in fact, it is non-negative),
$U_{\rm ad}\subset T^2(\Gamma_-)^3$ is convex and
$\lim_{\n{{\bf g}}_{T^2(\Gamma_-)^3}\to\infty,\ {\bf g}\in U_{\rm ad}}J({\bf g})=\infty$.
Therefore, there exists a unique minimum $\ol{{\bf g}}$ of $J$ in $U_{\rm ad}$
and the necessary and sufficient
conditions for the minimality of $\ol{{\bf g}}$ are (see e.g. \cite{lions71},  \cite{troltzsch})
\be\label{optim}
J'(\ol {\bf g})({w}-\ol {\bf g})\geq 0,\hspace{1.2cm} \forall {w}\in U_{\rm ad},
\ee
\be\label{optim-a}
\widetilde{\bf B}(\psi(\ol {\bf g}),v)=({\bf F}\ol {\bf g})(v),\quad \forall v\in \widehat{\mathbb{H}}.
\ee

Let
\[
(e_{\s T} h)(x)=\begin{cases}
h(x),\ &x\in {\s T}\cr 0,\ &x\in G\setminus {\s T}
\end{cases}
\]
for a function $h\in L^2({\s T})$
(that is, $e_{\s T} h$ is the extension by zero of $h$ onto $G$) and similarly we define extensions by zero for functions   $h\in L^2({\s C})$ and $h\in L^2({\s N})$ onto $G$.
The differential of $J'(\ol{\bf g})$ for $w\in T^2(\Gamma_-)^3$ is, recalling that ${\s D}\ol{\bf g}=D\psi(\ol{\bf g})$,
\bea\label{v10a:1}
J'(\ol{\bf g})w=&
2c_{\s T}\la D\psi(\ol{\bf g})-d_{\s T} ,D\psi(w)\ra_{L_2({\s T})}
+ 2c_{\s C}\la D\psi(\ol{\bf g})-d_{\s C} ,D\psi(w)\ra_{L_2({\s C})}\nonumber\\
&+ 2c_{\s N}\la D\psi(\ol{\bf g})-d_{\s N} ,D\psi(w)\ra_{L_2({\s N})}
+2c_{\rm sc}\la \ol{\bf g},w\ra_{T^2(\Gamma_-)^3}\nonumber\\
=&2c_{\s T}\la e_{\s T}(D\psi(\ol{\bf g})-d_{\s T}) ,D\psi(w)\ra_{L_2(G)}
+ 2c_{\s C}\la e_{\s C}(D\psi(\ol{\bf g})-d_{\s C}) ,D\psi(w)\ra_{L_2(G)}\nonumber\\
&+ 2c_{\s N}\la e_{\s N}(D\psi(\ol{\bf g})-d_{\s N}) ,D\psi(w)\ra_{L_2(G)}
+2c_{\rm sc}\la \ol{\bf g},w\ra_{T^2(\Gamma_-)^3}\nonumber\\
=&2c_{\s T}\la D^*e_{\s T}(D\psi(\ol{\bf g})-d_{\s T}) ,\psi(w)\ra_{L_2(G)^3}
+ 2c_{\s C}\la D^*e_{\s C}(D\psi(\ol{\bf g})-d_{\s C}) ,\psi(w)\ra_{L_2(G)^3}\nonumber\\
&+ 2c_{\s N}\la D^*e_{\s N}(D\psi(\ol{\bf g})-d_{\s N}) ,\psi(w)\ra_{L_2(G)^3}
+2c_{\rm sc}\la \ol{\bf g},w\ra_{T^2(\Gamma_-)^3}.
\eea

Denoting
\be\label{v12-a}
{\bf f}^*:=c_{\s T} D^*e_{\s T}(D\psi(\ol{\bf g})-d_{\s T}) 
+ c_{\s C} D^*e_{\s C}(D\psi(\ol{\bf g})-d_{\s C})
+ c_{\s N}D^*e_{\s N}(D\psi(\ol{\bf g})-d_{\s N} ),
\ee
we have ${\bf f}^*\in L^2(G\times S\times I)^3$, and so by Theorem \ref{ad-th} there
exists a unique $\psi^*\in \mathbb{H}$ such that
\be\label{v12}
\widetilde{\bf B}^*(\psi^*,v)=-\la {\bf f}^*,v\ra_{L^2(G\times S\times I)^3},\quad \forall v\in 
\widehat{\mathbb{H}}.
\ee

Let $w\in T^2(\Gamma_-)^3$. Then by the above notations
\be\label{v13}
\widetilde{\bf B}(\psi(w),v)=({\bf F}w)(v),\quad \forall v\in \widehat{\mathbb{H}}.
\ee
By Lemma \ref{le-sol-id} $\widetilde{\bf B}^*(\psi^*,\psi(w))$ and 
$\widetilde{\bf B}(\psi(w),\psi^*)$ are well-defined and (\ref{v12}) and (\ref{v13}) are valid with $v=\psi(w)$ and $v=\psi^*$, respectively that is,
\be\label{v13-a} 
\widetilde{\bf B}^*(\psi^*,\psi(w))=-\la {\bf f}^*,\psi(w)\ra_{L^2(G\times S\times I)^3}, 
\ee
\be\label{v13-b}
\widetilde{\bf B}(\psi(w),\psi^*)=({\bf F}w)(\psi^*)=\la \gamma_-(\psi^*), w\ra_{T^2(\Gamma_-)^3}.
\ee
Moreover, by (\ref{sol-id})
\be\label{v13-c}
\widetilde{\bf B}^*(\psi^*,\psi(w))=\widetilde{\bf B}(\psi(w),\psi^*).
\ee
Hence we have by  (\ref{v13-a}), (\ref{v13-b}), (\ref{v13-c})
\bea\label{v14}
&
2\la {\bf f}^*,\psi(w)\ra_{L^2(G\times S\times I)^3}=-2\widetilde{\bf B}^*(\psi^*,\psi(w))\nonumber\\
&
=-2\widetilde{\bf B}(\psi(w),\psi^*)=\la -2\gamma_-(\psi^*), w\ra_{T^2(\Gamma_-)^3}
\eea
Combining the previous expressions (\ref{v12-a}), \eqref{v10a:1} and \eqref{v14} thus leads to
\be\label{v14a}
J'(\ol{\bf g})w=&2\la {\bf f}^*,\psi(w)\ra_{L^2(G\times S\times I)^3}+2c_{\rm sc}\la \ol{\bf g},w\ra_{T^2(\Gamma_-)^3} 
=\la -2\gamma_-(\psi^*)+2c_{\rm sc}\ol{\bf g},w\ra_{T^2(\Gamma_-)^3}.
\ee

Note that $w+\ol{\bf g}\in U_{\rm ad}$ when $w\in U_{\rm ad}$.
Choosing for ${w}$ in the condition (\ref{optim})
subsequently $w+\ol{\bf g}$ and $0$, we see that
\bea
&J'(\ol {\bf g})w\geq 0,\quad \forall w\in U_{\rm ad} \label{v15:1}\\
&-J'(\ol {\bf g})\ol {\bf g}\geq 0.\label{v15:2}
\eea
From (\ref{v15:1}), (\ref{v15:2}) it follows that
\be\label{v15:3}
J'(\ol {\bf g})\ol {\bf g}= 0.
\ee
By (\ref{v14a}) and (\ref{v15:1}) one has
\be\label{v16}
\la-\gamma_-(\psi^*)+c_{\rm sc}\ol {\bf g},w\ra_{T^2(\Gamma_-)^3}\geq 0, \quad \forall w\in U_{\rm ad},
\ee
and so for each component $j=1,2,3$,
\be\label{v17}
-\gamma_-(\psi_j^*)+c_{\rm sc}\ol g_j\geq 0\ {\rm a.e.\ in}\ \Gamma_-.
\ee
Moreover, due to (\ref{v15:3}) we have
\be\label{v17a}
\la-\gamma_-(\psi^*)+c_{\rm sc}\ol {\bf g},\ol {\bf g}\ra_{T^2(\Gamma_-)^3}
=\sum_{j=1}^3\ol g_j(-\gamma_-(\psi_j^*)+ c_{\rm sc}\ol g_j)
=0,
\ee
and so by (\ref{v17}) for each $j=1,2,3$ (since also $\ol g_j\geq 0$),
\[
\ol g_j(-\gamma_-(\psi_j^*)+c\ol g_j)=0\ {\rm a.e.\ in}\ \Gamma_-.
\]
From this, using again (\ref{v17}),
one concludes that
\[
\ol g_j={1\over{c_{\rm sc}}}\max\{0,\gamma_-(\psi_j^*)\},\quad j=1,2,3,
\]
which is the claim (\ref{v9a}).
Finally, substituting this into the equation (\ref{optim-a}) we get from (\ref{optim-a})
 and (\ref{v12})
\[
\widetilde{\bf B}(\psi,v)&=({\bf F}\ol {\bf g})(v)=\la\ol {\bf g},\gamma_-(v)\ra_{T^2(\Gamma_-)^3},\\
\widetilde{\bf B}^*(\psi^*,v)&=-\la {\bf f}^*,v\ra_{L^2(G\times S\times I)^3}.
\]
Noticing that
\[
\la {\bf f}^*,v\ra_{L^2(G\times S\times I)^3}
=&c_{\s T} \la D\psi(\ol{\bf g})-d_{\s T},Dv\ra_{L^2({\s T})}
+ c_{\s C}\la D\psi(\ol{\bf g})-d_{\s C},Dv\ra_{L^2({\s C})} \\
&+ c_{\s N}\la D\psi(\ol{\bf g})-d_{\s N},Dv\ra_{L^2({\s N})},
\]
we get the system of equations (\ref{v10}) for the pair $(\psi, \psi^*)$. This completes the proof.
\end{proof}

\begin{remark}\label{re-init}
If we choose the whole space $\tilde{U}_{\rm ad}=T^2(\Gamma_-)^3$ as an admissible set
instead of $U_{\rm ad}$ (i.e. if non-negativity of admissible controls was not imposed),
we would find by considerations similar to those in the above proof
that the following variation of Theorem \ref{thinit} holds:

Under the assumptions of Theorem \ref{thinit},
the minimum $\min_{{\bf g}\in \tilde{U}_{\rm ad}}{J}({\bf g})$ exists and it is given by
\be\label{v9aa}
\ol{\bf g}={1\over{c_{\rm sc}}}\gamma_-(\psi^*)
\ee
where  $\psi^*\in \mathbb{H}$ is obtained from  the solution of the \emph{coupled linear system} of variational equations
\bea\label{v10a}
\widetilde{\bf B}^*(\psi^*,v)&+
c_{\s T}\la D\psi,Dv\ra_{ L^2({\s T})}+c_{\s C}\la D\psi,Dv\ra_{ L^2({\s C})}+c_{\s N}\la D\psi,Dv\ra_{ L^2({\s N})}\nonumber\\
&=c_{\s T}\la d_{\s T},Dv\ra_{ L^2({\s T})}+c_{\s C}\la d_{\s C},Dv\ra_{ L^2({\s C})}+c_{\s N}\la d_{\s N},Dv\ra_{ L^2({\s N})} \\
\widetilde{\bf B}(\psi,v)&=\la {1\over{c_{\rm sc}}}\gamma_-(\psi^*), \gamma_-(v)\ra_{T^2(\Gamma_-)^3},
\nonumber
\eea
for all  $v\in \widehat{\mathbb{H}}$.

By using this technique the initial solution for the global optimization
problem
would be taken to be
$\frac{1}{{c_{\rm sc}}}(\gamma_-(\psi^*))_+$. 
We point out that the equations (\ref{v10a}) are linear,
since $\psi^*\mapsto {1\over{{c_{\rm sc}}}}\gamma_-(\psi^*)$ is linear,
and therefore no iteration scheme is necessarily required in solving them.
Presumably, however, the solution of
non-linear optimization problem given in Theorem \ref{thinit}
should give a more accurate initial solution 
for the global optimization problem.
Similar observation is concerning for the internal therapy optimization described below.
 
\end{remark}

In the internal therapy ${\bf g}=0$ and so $\psi=\psi({\bf f}):=\psi({\bf f},0)$ and ${\s D}({\bf f}):=D(\psi({\bf f}))$. The object function for the initial solution may be
\be\label{va18}
J({\bf f})=&c_{\s T}\n{ d_{\s T}-{\s D}({\bf f})}_{L^2({\s T})}^2+c_{\s C}\n{ d_{\s C}-{\s D}({\bf f})}_{L^2({\s C})}^2 \\
&+c_{\s N}\n{d_{\s N}-{\s D}({\bf f})}_{L^2({\s N})}^2+c_{\rm sc}\n{{\bf f}}_{L^2(G\times S\times I)^3}^2. \nonumber
\ee 
and
\[
U'_{\rm ad}=\{{\bf f}\in L^2(G\times S\times I)^3\ |\ {\bf f}\geq 0\}.
\]
By arguments analogous to those used to prove Theorem \ref{thinit}
lead to the next result.

\begin{theorem}\label{thinit:internal}
Assume that the assumptions of Theorem \ref{thinit} hold.
Then the minimum $\min_{{\bf f}\in U'_{\rm ad}}{J}({\bf f})$ exists at the point ${\bf f}=\ol{{\bf f}}\in U'_{\rm ad}$
where
\be\label{v9a:3}
\ol{{\bf f}}=\frac{1}{{c_{\rm sc}}}(\psi^*)_+,
\ee
and  $\psi^*\in \mathbb{H}$
is obtained from the solution of the \emph{coupled non-linear} system of variational equations
\bea\label{v10:1}
\widetilde{\bf B}^*(\psi^*,v)&+
c_{\s T}\la D\psi,Dv\ra_{ L^2({\s T})}+c_{\s C}\la D\psi,Dv\ra_{ L^2({\s C})}
+c_{\s N}\la D\psi,Dv\ra_{ L^2({\s N})}\nonumber\\
&=c_{\s T}\la d_{\s T},Dv\ra_{L^2({\s T})}+c_{\s C}\la d_{\s C},Dv\ra_{L^2({\s C})}
+c_{\s N}\la d_{\s N},Dv\ra_{L^2({\s N})} \\
\widetilde{\bf B}(\psi,v)&=\la \frac{1}{{c_{\rm sc}}}(\psi^*)_+,v\ra_{ L^2(G\times S\times I)^3} \nonumber,
\eea
for all  $v\in \widehat{\mathbb{H}}$.
\end{theorem}

\begin{proof}
We shall content ourselves here with sketching briefly the part of proof leading to \eqref{v9a:3}.
Computations similar to those leading to \eqref{v16} in the proof of Theorem \ref{thinit},
would give in  the current context,
\bea\label{eq:optineq:internal}
\la -\psi^* + c_{\rm sc} \ol{\bf f},w\ra_{L^2(G\times S\times I)^3} \geq 0,\quad \forall w\in U'_{\rm ad},
\eea
hence
\[
-\psi^* + c_{\rm sc} \ol{\bf f}\geq 0\ \textrm{a.e. in}\ G\times S\times I,
\]
and those leading to \eqref{v17a} would give
\[
\la -\psi^* + c_{\rm sc} \ol{\bf f},\ol{\bf f}\ra_{L^2(G\times S\times I)^3}=0.
\]
Since $\ol{\bf f}\geq 0$ as $\ol{\bf f}\in U'_{\rm ad}$, 
we thus get
\bea\label{eq:opteq:internal}
\ol{\bf f}(-\psi^* + c_{\rm sc} \ol{\bf f})=0\ \textrm{a.e. in}\ G\times S\times I,
\eea
from which \eqref{v9a:3} easily follows.
\end{proof}

As we mentioned  these solutions can be utilized as an initial point for the global optimization schemes but they are not ready solutions for the treatment planning.
In \cite{tervo17-up}, at the end of section 8.2.3, we have discussed some other (special) admissible sets.

\begin{remark}
We finally remark that an optimization scheme can be formulated in such a way
that in external therapy the device (such as MLC) parameters, say ${\bf p}$, are {\it directly} as decision parameters both in static and dynamical  delivery techniques.
This is based on the fact that the incoming flux ${\bf g}$ can be expressed as a function of these parameters, say ${\bf g}={\bf g}({\bf p})$. 
Substitution of this expression ${\bf g}={\bf g}({\bf p})$ to 
${\s D}(0,{\bf g})$ gives the object function (\ref{io}) as a function of  ${\bf p}$. This approach would have obvious advantages, for example, device constraints could be taken into account directly in optimization.
The resulting object function is, however, highly multiextremal, and
thus seeking its global minimum is rendered more difficult.
\end{remark}

\section{Appendix}\label{appendix}

In this Appendix we expose (without proofs)
the change of variables from $(x,\omega,E)$-coordinates to $(x,v)$-coordinates where $v$ is the velocity of  particles.

\subsection{Transport problem in  space-velocity coordinates}\label{in-new-coord}

We perform  the change of variables where $(\omega,E)$ is replaced with the velocity variable $v$. 
Denote $r_0:=\sqrt{E_0},\ r_m:=\sqrt{E_m}$ and $B:=B(0,r_m)\setminus \ol B(0,r_0)$.
Let $S_{0}$ be the intersection of $S=S_{2}$ and the half plane ${\s T}_2:=\{x\in\R^3| x=(x_1,0,x_3),\ x_1\geq 0\}$  and 
let $B_0$ be the intersection of $B$ and ${\s T}_2$. Furthermore,
let $h:(S\setminus S_{0})\times I^\circ\to  B\setminus B_0$   be the diffeomorphism
\[
h(\omega,E):=\sqrt{E}\omega.
\]
This diffeomorphism (change of variables) is motivated by the dynamic equation $E={1\over 2}m |v|^2$ where $v$ is the velocity of particles. We find that
\[
h^{-1}(v)=\Big({v\over{|v|}},|v|^2\Big).
\]
Furthermore, define $H:G\times (S\setminus S_{0})\times I^\circ\to G\times (B\setminus B_0)$ by
\[
H(x,\omega,E):=(x,h(\omega,E))=(x,\sqrt{E}\omega)=(x,v).
\]
Then $H$ is a ($C^\infty$) diffeomorphism and 
\[
H^{-1}(x,v)=(x,{v\over{|v|}},|v|^2).
\]

Let
\be\label{orib12a}
T\psi:=a(x,E){\p \psi{E}}
+b(x,E)\Delta_S\psi
+\omega\cdot\nabla_x\psi+\Sigma\psi-K_{r}\psi.
\ee
We express  the transport problem 
\be\label{b-in-cond}
T\psi=f \ {\rm in}\ G\times S\times I,\ \psi_{|\Gamma_-}=g(y,\omega,E),\ \psi(.,.,E_m)=0
\ee
in  $v$-coordinates.
The equation (\ref{orib12a}) is equivalent to
\bea\label{orib13a}
&
\Big(a{\p {\psi}E}\Big)\circ H^{-1}+\big(b(x,E)\Delta_S\psi\big)\circ H^{-1}+
(\omega\cdot\nabla_x\psi)\circ H^{-1}\nonumber\\
&
+(\Sigma\psi)\circ H^{-1}-(K_r\psi)\circ H^{-1}=f\circ H^{-1} 
\eea
in $G\times (B\setminus B_0)$.
The terms  
$\Big(a{\p {\psi}E}\Big)\circ H^{-1}$, $\big(b(x,E)\Delta_S\psi\big)\circ H^{-1}$,
$(\omega\cdot\nabla_x\psi)\circ H^{-1}$,
$(\Sigma\psi)\circ H^{-1}$, $(K_r\psi)\circ H^{-1}$ can be computed where the most work is required for the term $\big(b(x,E)\Delta_S\psi\big)\circ H^{-1}$.  

Let $\Psi:=\psi\circ H^{-1}$ or equivalently $\psi=\Psi\circ H$.
The resulting transport equation in the new coordinates is given by
\bea\label{f-app-op-v}
&
T\Psi=
\tilde b(x,v)P(x,v,D)\Psi+
(\tilde a(x,v)-2\tilde b(x,v))(v\cdot \nabla_v\Psi)\nonumber\\
&
+{1\over{|v|}}(v\cdot\nabla_x\Psi)+\widetilde\Sigma\Psi-\widetilde K_{r}\Psi
\eea
where
\be\label{ta}
\tilde a(x,v)={1\over{2|v|^2}}a(x,|v|^2),\
\tilde b(x,v)=b(x,|v|^2),
\ee
\be\label{ts}
\widetilde\Sigma(x,v):=(\Sigma\circ H^{-1})(x,v)
=\Sigma(x,|v|^2)
\ee
and  $P(x,v,D)\Psi$ is (here $\partial_{ij}\Psi={{\partial^2\Psi}\over{\partial v_i\partial v_j}}$)
\bea\label{P}
&
(P(x,v,D)\Psi)(x,v)\nonumber\\
&
:=
{{|v|^2v_2^2+v_1^2v_3^2}\over{{v_1^2}+v_2^2}}\ 
(\partial_{11}\Psi)(x,v)-2v_1v_2(\partial_{12}\Psi)(x,v)
-2v_1v_3(\partial_{13}\Psi)(x,v)
-2v_2v_3(\partial_{23}\Psi)(x,v)
\nonumber\\
&
+
{{|v|^2v_1^2+v_2^2v_3^2}\over{{v_1^2}+v_2^2}}\
(\partial_{22}\Psi)(x,v)
+(v_1^2+v_2^2)(\partial_{33}\Psi)(x,v)
. 
\eea
The restricted collision operator is
\[
\widetilde K_r\Psi=\widetilde K_r^1\Psi+\widetilde K_r^2\Psi+\widetilde K_r^3\Psi
\]  
where
\be\label{co-4}
(\widetilde K_r^1\Psi)(x,v):=
((K_r^1\psi)\circ H^{-1})(x,v)
=
\int_{B'}
\tilde\sigma^1(x,v',v)
\Psi(x,v') dv' ,
\ee
\be\label{co-5}
(\widetilde K_r^2\Psi)(x,v):=
((K_r^2\psi)\circ H^{-1})(x,v)
=
\int_{B'}
\tilde\sigma^2(x,v',v)
\Psi(x,v') dv' ,
\ee
\be\label{co-3}
(\widetilde K_r^3\Psi)(x,v):=
((K_r^3\psi)\circ H^{-1})(x,v)
=
\int_{B'}\int_{0}^{2\pi}
\tilde\sigma^3(x,v',v)
\Psi(x,\tilde\gamma(v',v)(s))ds dv' .
\ee
In the above formulas the restricted differential cross-sections are
\[
\tilde\sigma^1(x,v',v):={2\over{|v'|}}\sigma^1(x,{{v'}\over{|v'|}},{{v}\over{|v|}}
,|v'|^2,|v|^2),
\]
\[
\tilde\sigma^2(x,v',v):={2\over{|I'|}}{1\over{|v'|}}\sigma^2(x,{{v'}\over{|v'|}},{{v}\over{|v|}},|v|^2)
\]
and 
\[
\tilde\sigma^3(x,v',v):={1\over{2\pi}}{1\over{|v'|}}\sigma^3(x,|v'|^2,|v|^2).
\]
In addition,
\[
\tilde\gamma(v',v)(s):=|v'|\gamma(|v'|^2,|v|^2,{v\over{|v|}})(s).
\]

Let
\[
\widetilde \Gamma_-:=\{(y,v)\in \partial G\times B|\ v\cdot\nu(y)<0\}.
\]
Then we see that the inflow boundary condition for $\Psi$ is
\be\label{fpr-1}
\Psi_{|\widetilde \Gamma_-}=\tilde g
\ee
where $\tilde g(y,v):=g(y,{v\over{|v|}},|v|^2)$.
The initial condition for $\Psi$ is
\be\label{fpr-2}
\Psi(.,v)=0\ {\rm for}\ |v|=r_m 
\ee 
where $r_m=\sqrt{E_m}$.

Via the use of these $(x,v)$-coordinates one can avoid the calculations on the  (compact)  manifold $S$ (which may be challenging). 
With the same kind of analysis as above one can show
the existence result of solutions (analogous to Theorem \ref{csdath3}) for the problem 
\[
T\Psi=\tilde f,\ \Psi_{|\widetilde \Gamma_-}=\tilde g,\ 
\Psi(.,v)=0\ {\rm for}\ |v|=r_m .
\]
Similarly the solution for the corresponding inverse problem (analogous to Theorem \ref{thinit}) can be verified. We omit all formulations here.

\end{document}